\documentclass [12pt] {report}
\usepackage[dvips]{graphicx}
\usepackage{epsfig, latexsym,graphicx}
\usepackage{amssymb}
\newcommand {\D}[2] {\displaystyle\frac{\partial{#1}}{\partial{#2}}}

\newcommand {\ga} {\gamma}

\newcommand {\si} {\sigma}
\newcommand {\Si} {\Sigma}
\newcommand {\de} {\delta}

\newcommand {\fr} {\displaystyle\frac}
\newcommand {\wt} {\widetilde}
\newcommand {\be} {\begin{equation}}
\newcommand {\ee} {\end{equation}}
\newcommand {\ba} {\begin{array}}
\newcommand {\ea} {\end{array}}
\newcommand {\bp} {\begin{picture}}
\newcommand {\ep} {\end{picture}}
\newcommand {\bc} {\begin{center}}
\newcommand {\ec} {\end{center}}
\newcommand {\bt} {\begin{tabular}}
\newcommand {\et} {\end{tabular}}
\newcommand {\lf} {\left}
\newcommand {\rg} {\right}

\newcommand {\cF} {{\cal F}}

\newcommand {\cR} {{\cal R}}
\newcommand {\cS} {{\cal S}}

\newcommand {\ses} {\medskip}

\newcommand {\e} {\mathop{\rm e}\nolimits}

\newcommand {\bibit} {\bibitem}
\newcommand {\nin} {\noindent}

\newcommand {\Rho} {\mbox{\large$\rho$}}

\newcommand {\sign} {\mathop{\rm sign}\nolimits}

\newcommand {\wh} {\widehat}

\newcommand {\cA} {{\cal A}}

\usepackage{epsfig, latexsym,graphicx}
 \usepackage{amsmath}

\newcommand {\cD} {{\cal D}}

\def\2#1#2#3{{#1}_{#2}\hspace{0pt}^{#3}}
\def\3#1#2#3#4{{#1}_{#2}\hspace{0pt}^{#3}\hspace{0pt}_{#4}}
\newcounter{sctn}
\def\sec#1.#2\par{\setcounter{sctn}{#1}\setcounter{equation}{0}
                  \noindent{\bf\boldmath#1.#2}\bigskip\par}

\begin {document}

\begin {titlepage}

\vspace{0.1in}

\begin{center}

{\Large \bf  Finslerian angle-preserving  connection
}
\ses
{\Large \bf in two-dimensional  case.    Regular realization}

\end{center}

\vspace{0.3in}

\begin{center}

\vspace{.15in} {\large G.S. Asanov\\} \vspace{.25in}
{\it Division of Theoretical Physics, Moscow State University\\
119992 Moscow, Russia\\
{\rm (}e-mail: asanov@newmail.ru{\rm )}} \vspace{.05in}

\end{center}

\begin{abstract}

\ses

We show that the metrical connection can be introduced in the
two-dimensional Finsler space such that entailed parallel transports
along curves joining points of the underlying manifold
keep the two-vector angle as well as the length of the tangent vector,
thereby realizing isometries of tangent spaces under the
parallel transports.
The curvature tensor is found.
In case of the Finsleroid-regular space,
constructions possess  the $C^{\infty}$-regular status globally
regarding
the dependence on tangent vectors.
Many involved and important relations are explicitly derived.

\medskip


\end{abstract}

\end{titlepage}

\vskip 1cm

\ses

\ses

\setcounter{sctn}{1} \setcounter{equation}{0}

\bc  {\bf 1. Motivation and description}  \ec

\bigskip

During all the history of development of the Finsler geometry
the notion of connection
was attracted  sincere and great attention of investigators devoted to general theory
as well as  to specialized applications.
The  methods of construction of connection are founded upon
setting forth a  convenient system of axioms.
Various standpoints were taken to get deeper insights into the notion
(see  [1-5]).

The general idea underlining the   present work is to
 set forth
 the requirement that the connection  be compatible with the preservation
of the two-vector angle under the  parallel transports of vectors.

The notion of  angle is of key significance in geometry.
In the field of two-dimensional Finsler spaces
the  angle between two vectors of a given tangent space
can naturally be measured by the area of the domain
bounded by  the vectors and the indicatrix arc.
The theorem can be proved which states
that a diffeomorphism between two Finsler spaces
is an isometry iff it keeps the  angles thus appeared.
This fundamental Tam{\' a}ssy's theorem [6],
which  explains  us that
the angle structure fixes the metric structure in the
 Finsler space,
  gives rise to the following  important question:
{\it Does the angle structure also generate  the connection?}
The present work proposes a positive and
explicit answer,
 confining the Finslerian consideration to the two-dimensional case.

{%\pgbrk}

Let $M$ be a
$C^{\infty}$-differentiable   2-dimensional
manifold, $ T_xM$ denote the tangent space  supported by the  point $x\in M$,
and $y\in T_xM\backslash 0$  mean tangent vectors.
Given a Finsler metric function $F=F(x,y)$, we obtain the  two-dimensional
 Finsler space
$\cF_2 = (M, F)$.

We shall use the standard Finslerian notation for local components
 $l^k=y^k/F, ~ y_k=F\partial F/\partial y^k\equiv g_{kn}y^n,~
g_{ij}=\partial y_i/\partial y^j$ of the unit vector,
 the covariant tangent vector, and the Finsler metric tensor, respectively.
The covariant components
$l_k=g_{kn}l^n$
can be obtained from
$l_k=\partial F/\partial y^k$.
By means of the contravariant components
$g^{ij}$
the Cartan tensor $C_{ijk}=(1/2)\partial g_{ij}/\partial y^k$
can be contracted to yield
the vector $C_k=g^{ij}C_{ijk}$.
It is convenient to use the tensor
$A_{ijk}=FC_{ijk}$
and the vector $A_k=FC_k=g^{ij}A_{ijk}$.
The indices $i,j,...$ are specified over the range (1,2).
The square root
$\sqrt{}$
stands always in the positive sense.
It is often convenient to apply the expansion
$A_{ijk}=Im_im_jm_k$ in terms of the $m_i$ obtainable
from $g_{ij}=l_il_j+m_im_j$, where $I$ thus appeared is the so-called
main scalar.
Our consideration will be of local nature, unless otherwise is
stated explicitly.

\ses

 {%\pgbrk}

To each point  $x\in M$,
 the Finsler  space  $\cF_2$ associates
the tangent Riemannian  space, to be denoted by
$\cR_{\{x\}}~:=\{T_xM, g_{ij}(x,y)\}$, in which
  $x$ is treated fixed and  $y\in T_xM$ is variable.
In  the Riemannian space the  $\cR_{\{x\}}$
reduces  to
the tangent Euclidean space.
The remarkable and well-known property of the
Riemannian Levi Civita connection is that the entailed parallel transports
along curves drawn on
the underlined  manifold
 keep the length of the tangent vectors
and produce the isometric mapping of
the tangent Euclidean spaces.

We show that these two fundamental Riemannian properties
can successfully be extended to operate in  the Finsler space  $\cF_2$.
Namely, if sufficient smoothness holds then
 it proves possible to introduce  the
respective  connection coefficients
 $\{N^k_i(x,y), D^k{}_{in}(x,y)\}$
 in a simple and  explicit way.
The coefficients
$N^k_i(x,y)$ are required  to construct the conventional operator
$
d_n=\partial_{x^n}+N^k_n(x,y)\partial_{y^k},
$
where
 $\partial_{x^n}=\partial/\partial{x^n}$
and
$\partial_{y^m}=\partial/\partial{y^m}$.
The keeping
of
 the Finsler    length  of the tangent vectors
 means $d_nF=0$.
Let us attract also the angle function $\theta=\theta(x,y)$
(to  measure  the length $ds$ of  infinitesimal arc on the indicatrix
by $d\theta$)
and raise forth the requirement that
$d_n\theta=k_n$ with  a  covariant vector field
$k_n=k_n(x)$.
If this is  fulfilled,
then for pairs $\{\theta_1=\theta(x,y_1),\theta_2=\theta(x,y_2)\}$
we obtain the nullification $d_n(\theta_2-\theta_1)=0$
which tells us that
the  {\it   preservation of the two-vector angle}
$\theta_2-\theta_1$ holds true under the parallel transports initiated by
the coefficients
$N^k_i(x,y)$.

There arise   the  coefficients
$D^k{}_{in}(x,y)=-N^k{}_{in}(x,y)$
 with $N^k{}_{in}(x,y)=\partial N^k_i(x,y)/\partial y^n$.
 A careful analysis
has shown that a simple and attractive proposal of
the  coefficients $N^k_i(x,y)$
(namely, (2.6) of Section 2)
can be made such that
nullifications
$d_iF=0$
and
$d_i\theta-k_i=0$
are simultaneously  satisfied,
and also
the vanishing
$l_kN^k{}_{nmi}=0$
holds fine because of the
representation
$
FN^k{}_{nmi}=
-A^k{}_{mi}d_n\ln|I|
$
entailed (see (2.14)),
where
  $N^k{}_{nmi}=\partial N^k{}_{nm}/\partial y^i$,
  so that the action of the arisen covariant derivative
on the involved Finsler metric tensor
yields just the zero.
The coefficients $D^k{}_{in}$
are not symmetric in the subscripts $i,n$.

Having realized this program (in Section 2),
 we feel sure that the arisen mappings of the space
 $\cR_{\{x\}}$
 under the respective parallel transports along the curves running on  $M$
are isometries.

The coefficients  $\{N^k_i(x,y), D^k{}_{in}(x,y)\}$ obtained in this way
are {\it not} constructed from the Finsler
metric tensor and derivatives of the tensor.
This circumstance may be estimated to be a cardinal  distinction
of the Finsler connection induced by  the angle structure from the
conventional Riemannian precursor which exploits
the Riemannian Christoffel symbols to be
 the coefficients
$D^k{}_{in}$.
The structure of the coefficients $N^k_n$ involves the derivative
$\partial \theta/\partial x^n$ on the equal footing with the derivative
$\partial F/\partial x^n$
(see (2.6) in Section 2).

\ses

{%\pgbrk}

The involved vector field $k_i=k_i(x)$ may be taken arbitrary.
However, the field can be specified
if the Riemannian limit of the connection proposed is attentively considered.
Indeed,
in the Riemannian limit
 the connection coefficients   $D^k{}_{nh}(x,y)$
reduce to the coefficients
  $\bar L^k{}_{nh}=-L^k{}_{nh}=\bar L^k{}_{nh}(x)$ which are
{\it not symmetric} with respect to the
subscripts (see (4.6)).
 If we want to obtain the torsionless  coefficients,
like to the Riemannian geometry proper,
 we must
 make the choice
$k_n=-n^h\nabla_n  {\wt b}_h$
in accordance with (4.16),
where
$\nabla_n $ is the Riemannian covariant derivative
taken with
the Christoffel symbols  $a^k{}_{nh}$.
The
$\wt b_h=\wt b_h(x)$ is a vector field chosen to fulfill
$\theta(x,\wt b(x))=0$,
and the pair $n_i,\wt b_i$ is orthonormal.
With the choice we obtain
 $\bar L^k{}_{nh}= a^k{}_{nh}$,
thereby completely specifying
the coefficients
 $\{N^k_i, D^k{}_{in}\}$.

 It is appropriate to construct the
{\it osculating } Riemannian metric tensor
along the vector field
${\wt b}^i={\wt b}^i(x)$
and introduce
the $\theta$-{\it associated Riemannian space}
to  compare the Finsler properties of the space
$\cF_2$ with the properties of the Riemannian precursor.

 {%\pgbrk}

\ses

\ses

In Section 2
the required coefficients $N^k_n$ are proposed
and nearest implications are indicated.
By the help of the identities
$
\partial m^k/\partial y^m=  -(Im^k+l^k)m_m/F$ and
$\partial m_m/ \partial y^i= (Im_m-l_m)m_i/F$,
the validity of   the vanishing
$l_kN^k{}_{nmi}=0$
can readily be verified.
The angle function $\theta$ introduced does measure
the length of the indicatrix arc
according to
$ds=d\theta$ (see (2.23)).
We derive also the equality
$
\Si_{\{y_1,y_2\}}=(1/2)(\theta_2-\theta_1),
$
where the left-hand side is the area of the sector
bounded by  the vectors
$y_1,y_2$ and the indicatrix arc  (see (2.27)).
The equality
demonstrates clearly
 that,
 in context of the two-dimensional  theory to which our treatment is restricted,
  the  method of introduction of  angle
by the help of the function $\theta$
is equivalent to the method founded in [6] on the notion of area.
We are entitled therefore to raise the thesis  that,
in such a context,
the angle-preserving connection is tantamount to the
area-preserving connection.

In Section 3
we show
how
the curvature tensor $\Rho_k{}^n{}_{ij}$ of the space $\cF_2$
can be explicated from the commutator of the covariant derivative arisen,
yielding  the astonishingly  simple representation
(3.14).

\ses

In Section 4 we outline   Riemannian counterparts.
It appears that the tensor $\Rho_{knij}=g_{nm}\Rho_k{}^m{}_{ij}$
is factorable,
in accordance with (4.19).

\ses

In Conclusions
several ideas are emphasized.

\ses

In Appendix A the validity of
succession of  important relations appeared in the analysis
has been explicitly demonstrated.

{%\pgbrk}

The possibility of global realization of
the angle-preserving  connection
implies high regularity properties of the Finsler metric function
and the angle function.
Such a lucky possibility  occurs
in the Finsleroid-regular space
$\cF\cR^{PD}_{g;c}$
(introduced and studied in [7-10]).
 Appendix B is devoted to  the space.
The Finsler metric function
$K=K(x,y)$ of the  space
$\cF\cR^{PD}_{g;c}$
 involves a Riemannian metric tensor
 $a_{mn}$ and
the vector field $b^n=b^n(x)$  which
 represents the distribution of axis of  indicatrix.
 We have two scalars,
 namely
the characteristic scalar $g=g(x)$ and
the norm
$c=c(x)=||b||\equiv\sqrt{a^{mn}b_mb_n}$
 of the 1-form $b=b_i(x)y^i$.
The metric function $K$ is not absolute homogeneous.

Simple direct calculation shows that the partial derivative
$\partial K/\partial x^n$
obeys the   total regularity with respect to the vector variable $y$
(see (B.59)).
The same regularity is shown by
the partial derivative
$\partial \theta/\partial x^n$
of the involved
angle function $\theta=\theta(x,y)$
(see (B.83)).
Therefore, all the ingredients in  the coefficients
$N^k_n$ of the form
 proposed by (2.6)
 are of this high regularity.
Thus we observe the remarkable  phenomenon that the
 space
$\cF\cR^{PD}_{g;c}$
possesses the  angle-preserving connection
of the $C^{\infty}$-regular status  globally regarding the $y$-dependence.
All the formulas appeared in Appendix B are such that their
right-hand parts are of this high  regularity
status.
  Arbitrary (smooth) dependence on $x$ in  $g=g(x)$,  $b_i=b_i(x)$,
  and
  $a_{ij}=a_{ij}(x)$
   is permitted.
The particular condition $c=const$ simplifies many representations.
With  $c=const$ we make an attractive proposal (B.85)
for the  coefficients
$N^k_n$ and demonstrate  in the step-by-step way
that  all the three conditions listed in (2.8)
are obtained and the  representation
of the announced form (2.6) is fulfilled.
The curvature tensor
has been evaluated, yielding the representations of the forms
(3.15) and  (4.19)
with the global  $C^{\infty}$-regular status of dependence on $y$.
Using an appropriate regular atlas of charts in the space
$\cR_{\{x\}}~:=\{T_xM, g_{ij}(x,y)\}$,
it proves possible to
verify that over  all the space $T_xM\backslash0$
the function
$\theta=\theta(x,y)$
is  smooth
  of class $C^{\infty}$ with respect to $y$.
The entailed two-vector angle
$\theta_2-\theta_1$ is symmetric and additive.
 The $\theta$  is represented by means  of integral and is not obtainable
through composition of elementary functions.
Several representations essential for evaluations
are indicated.

Quite similar evaluation can be performed for the Randers metric function
(Appendix C), yielding again
the  angle-preserving connection
of the $C^{\infty}$-regular status  globally regarding the $y$-dependence.

{%\pgbrk}

\setcounter{sctn}{2} \setcounter{equation}{0}

\ses

\ses

\bc
{\bf 2. Proposal of connection coefficients}
\ec

\ses

\ses

It is convenient to proceed with the skew-symmetric tensorial object
$\epsilon_{ik}=\sqrt{\det(g_{mn})}\,\ga_{ik}$, where
$\ga_{11}=\ga_{22}=0$ and $\ga_{12}=-\ga_{21}=1,$
to construct
\be
m_i=
-
\epsilon_{ik}l^k.
\ee
The angular metric tensor $h_{ij}=g_{ij}-l_il_j$
and the Cartan tensor $C_{ijk}$
are factorized,
and
the Finsler  metric tensor $g_{ij}$ is expanded,
 according to
\be
h_{ij}=m_im_j, \qquad  A_{ijk}=Im_im_jm_k, \qquad
g_{ij}=l_il_j+m_im_j.
\ee
It is also convenient
to  introduce the  $\theta=\theta(x,y)$ by the help of the equation
\be
F\D{\theta}{y^n}=m_n,
\ee
assuming that the  function $\theta$
is positively homogeneous of degree zero with respect to
the variable $y$.
These formulas are known from Section 6.6  of the book [1].
We denote the main scalar by $I$, instead of  $J$ used in the book.
Our $\theta$ is the $\varphi$ of  Section 6.6  of  [1].
The object
$\epsilon_{ik}$ is a pseudo-tensor, whence $m_i$ is a pseudo-vector
and
$I,\theta$
are pseudo-scalars.
However, we don't consider the coordinate reflections
and, therefore, we are entitled to refer to these objects as to
``the vector $m_i$" and ``the scalars $I,\theta$".
{%\pgbrk}

We need  the coefficients
$N^k_n=N^k_n(x,y)$
  to construct the operator
\be
d_n=
\D{}{x^n}+N^k_n   \D{}{y^k}
\ee
which generates a covariant vector $d_n W$
when is applied  to an arbitrary
 differentiable  scalar $W=W(x,y)$.
We shall use also   the derivative coefficients
\be
N^k{}_{nm}=\D{N^k_{n}}{y^m}, \qquad  N^k{}_{nmi}=\D{N^k{}_{nm}}{y^i}.
\ee

\ses

\ses

PROPOSAL.
{\it Take   the coefficients $N^k_n$ according to
the expansion
\be
N^k_n
=
-l^k    \D{F}{x^n}
-
F  m^k
\breve P_n
 \ee
with
\be
\breve P_n= \D{\theta}{x^n}- k_n,
\ee
\ses
where
$k_n=k_n(x)$ is a covariant vector field,
such that  the equalities
\be
d_nF=0, \qquad
d_n\theta=k_n, \qquad
y_kN^k_{nmi}=0
\ee
be realized.
}

\ses

\ses

The vanishing $d_nF=0$ and the equality
$ d_n\theta=k_n$
just follow from the choice (2.6).
Considering two values
  $\theta_1=\theta(x,y_1)$  and
 $\theta_2=\theta(x,y_2)$,
we  have
 \ses
\be
d_n\theta_1=\D{\theta_1}{x^n}+N^k_n(x,y_1)  \D{\theta_1}{y_1^k},
\qquad
d_n\theta_2=\D{\theta_2}{x^n}+N^k_n(x,y_2)
  \D{\theta_2}{y_2^k},
\ee
and from
$ d_n\theta=k_n$
 we may
conclude
that
 {\it the}
 {% angle}
{\it  preservation}
\be
d_n(\theta_2-\theta_1)=0
\ee
{\it holds because  the vector field $k_n$ is independent of tangent vectors $y$}.

From (2.6) it follows directly that
\be
N^k{}_{nm}
=
-
l^k \D{l_m}{x^n}
-
l_m  m^k
\breve P_n
-
F \D{ m^k}{y^m}
\breve P_n
-
 m^k
\D{m_m}{x^n}.
\ee

{%\pgbrk}

It is convenient to use the identities
\be
F\D{m_k}{y^m}=-l_km_m+Im_mm_k, \qquad
F\D{m^k}{y^m}=-Im^km_m-l^km_m
\ee
(they are tantamount to the identities written
 in  formula (6.22) of chapter 6 in [1]),
together with their immediate implication
$$
F\D{(m^km_n)}{y^m}=-(l_nm^km_m+l^km_nm_m).
$$
\ses
Short evaluations show that
\be
N^k{}_{nmi}=\fr1Fm^km_m \lf(F\D I{y^i}
\breve P_n
-
m_i
\D I{x^n}
\rg).
\ee
 Indeed,
from (2.11) it follows straightforwardly that
$$
N^k{}_{nmi}=
-\fr1F h^k_i\partial_nl_m
-l^k\partial_n\lf(\fr1Fh_{mi}\rg)
-\fr1F h_{mi}m^k\breve P_n
+l_m
 \Bigl( {I}\,m^km_i+l^km_i
\Bigr)
\breve P_n
-l_mm^k\partial_n\lf(\fr1Fm_i\rg)
$$

\ses

$$
+  \D{{I}}{y^i}\,m^km_m \breve P_n
+ {I}\,m^km_m
\partial_n\lf(\fr1Fm_i\rg)
-
Il^km_i
m_m \breve P_n
-  Im^kl_mm_i
\breve P_n
$$

\ses

\ses

$$
+\fr1F h^k_im_m \breve P_n
-l^k
(l_m-Im_m)m_i
\breve P_n
+l^km_m\partial_n\lf(\fr1Fm_i\rg)
$$

\ses

\ses

$$
+( Im^k+l^k)m_i
\fr1F
\partial_nm_m
+m^k\partial_n
\Bigl(\fr1Fl_mm_i-\fr1F{I}\,m_im_m\Bigr),
$$
where $\partial_n$ means $\partial/\partial x^n$.
The identity $h_{mi}=m_mm_i$ has been taken into account.
Canceling similar terms leads to (2.13).

Owing to  the identity $l_km^k=0$,
the vanishing
$l_kN^k{}_{nmi}=0$ holds fine.

\ses

Because of $y^i\partial_{y^i}I=0$,
the equality (2.13) can be written in the concise form
\be
N^k{}_{nmi}=-\fr1Fm^km_mm_i
d_nI
\equiv
-\fr1FA^k{}_{mi}d_n\ln|I|.
\ee

{%\pgbrk}

The sought {\it Finsler connection}
\be
{\cal FC}=\{N^k_n,D^k{}_{nm}\}
\ee
involves also the coefficients
$D^k{}_{nm}=D^k{}_{nm}(x,y)$
which are  required to construct the operator
of   the  covariant derivative $\cD_n$ which action
is exemplified in  the conventional way:
\be
\cD_nw^k{}_m~:=   d_nw^k{}_m  +   D^k{}_{nh}w^h{}_m   -  D^h{}_{nm}w^k{}_h,
\ee
where
$w^k{}_m=w^k{}_m(x,y)$ is an arbitrary differentiable (1,1)-type tensor.

If we differentiate the vanishing $d_iF=0$
with respect to the variable $y^j$ and
multiply the result by $F$, we obtain the vanishing
\be
\cD_iy_j~:=\D{y_j}{x^i}+N_i^kg_{kj} - D^h{}_{ij}y_h=0
\ee
\ses
when the choice
\be
D^k{}_{in} =-N^k{}_{in}
\ee
is made. Differentiating (2.17) with respect to  $y^n$
just manifests that the choice is also of success
 to fulfill  the metricity condition
\be
\cD_ig_{jn}~:=
d_ig_{jn} - D^h{}_{ij}g_{hn}- D^h{}_{in}g_{jh}
=0,
\ee
because of $y_kN^k_{nmi}=0$.

If we contract (2.18) by $y^n$ and take into account the definition of the coefficients
$N^k{}_{in}$
indicated in (2.5),
 we obtain the equality
\be
N^k_i=-D^k{}_{in}y^n.
\ee

\ses

The contravariant version of the vanishing (2.17) is obtained through the chain
\be
\cD_iy^j~:=d_iy^j +D^j{}_{ih}y^h=
N_i^j +D^j{}_{ih}y^h=0.
\ee

Because of
$\cD_ih^{nk}=0$, applying the derivative $\cD_i$ to
the equality
$h^{nk}=m^nm^k$  (see (2.2))
 and contracting the result by $m_n$, we conclude that
\be
\cD_im^k=0, \quad \text{which means} ~ ~ d_im^k=N^k{}_{ih}m^h.
\ee

{%\pgbrk}

 \ses

Because of the homogeneity, the unit tangent vector components $l^n=l^n(x,y)$ can
obviously be regarded as  functions $l^n=L^n(x,\theta)$
of the pair $(x,\theta)$.
Let us denote $l^n_{\theta}=\partial L^n/\partial \theta$.
Since
$\partial F/\partial \theta=0$
and
$l_nl^n_{\theta}=0$,
we may conclude from (2.3) that  $l^n_{\theta}=m^n$.
Measuring the length of the indicatrix (which is defined by $F=1$)
by means of a parameter  $s$, so that
$ds=\sqrt{g_{ij}dl^idl^j}$,
we obtain
$ds=\sqrt{g_{ij}l^i_{\theta} l^j_{\theta}}\,d\theta=d\theta$,
assuming $ds>0$ and  $d\theta>0$.
Thus
\be
ds=d\theta \quad \text {along the indicatrix},
\ee
which explains us that the $\theta$ {\it  measures the length of indicatrix.}

If at a fixed $x$ we introduce in the tangent space $T_xM$ the coordinates
$z^A=\{z^1=F, z^2=\theta\}$ and consider the respective transforms
\be
G_{AB}=g_{ij}\D{y^i}{z^A}\D{y^i}{z^B}, \quad A=1,2, ~ B=1,2,
\ee
of the Finsler metric tensor components $g_{ij}$,
we obtain simply
\be
G_{11}=1, \quad G_{12}=0, \quad G_{22}=F^2.
\ee
With these components, it is easy to calculate the area of domain of the
tangent space  $T_xM$
by using the integral
 measure
 \be
\int\sqrt{\det(G_{AB})}\,dz^1dz^2=\int FdFd\theta.
\ee
\ses
In particular, for the sector
$\si_{\{y_1,y_2\}}\subset T_xM$ bounded by the vectors
$y_1,y_2$ and the indicatrix arc
we obtain by integration
the area $\Si_{\{y_1,y_2\}}$ which is given by
\be
\Si_{\{y_1,y_2\}}=\fr12(\theta_2-\theta_1),
\ee
so that in the two-dimensional case the angle in the Finsler geometry can be
 defined by the area just in the same way as
in the Riemannian geometry.
The difference
\be
\theta_2-\theta_1=\theta(x,y_2)-\theta(x,y_1)
\ee
can naturally be regarded as the value of angle between two vectors
$y_1,y_2\in T_xM$, the {\it two-vector angle} for short.
The formula (2.27) is equivalent to Definition (2) of  [6] which was proposed
to define angle by area; we use the right orientation of angle.

\ses

As well as the area is attributed to the tangent space by means of the
integral measure
(2.26)
and the conditions
$
d_nF=0$ and $d_n(\theta_2-\theta_1)=0
$
are fulfilled,
{\it  the angle-preserving connection keeps the area
under parallel transports}
 along curves joining point to point
in the background  manifold.
Thus we are entitled to set forth the thesis:
{\it  the angle-preserving connection  is  the area-preserving connection}.

{%\pgbrk}

\ses

\setcounter{sctn}{3} \setcounter{equation}{0}

\ses

\ses

\bc{
\bf 3.  Curvature tensor}
\ec

\ses

\ses

With
arbitrary coefficients
$\{N^k_n,D^k{}_{nm}\}$,
commuting the covariant derivative (2.16) yields the equality
\be
\lf(\cD_i\cD_j-\cD_j\cD_i\rg) w^n{}_k
=M^h{}_{ij}\D{w^n{}_k}{y^h}  -E_k{}^h{}_{ij}w^n{}_h
+E_h{}^n{}_{ij}w^h{}_k
\ee
with the   tensors
\be
 M^n{}_{ij}~:= d_iN^n_j-d_jN^n_i
\ee
and
\be
E_k{}^n{}_{ij}~: =
d_iD^n{}_{jk}-d_jD^n{}_{ik}+D^m{}_{jk}D^n{}_{im}-D^m{}_{ik}D^n{}_{jm}.
\ee
\ses

If the choice
$ D^k{}_{in} =-N^k{}_{in} $
is made (see (2.18)),
 the tensor (3.2) can be written in the form
\be
 M^n{}_{ij}= \D{N^n_j}{x^i}-\D{N^n_i}{x^j}
-N_i^hD^n{}_{jh}+N_j^hD^n{}_{ih},
\ee
which entails
the equality
\be
E_k{}^n{}_{ij}= -\D{M^n{}_{ij}}{y^k}.
\ee

By applying the commutation rule (3.1)
to the vanishing set
$\{\cD_iF=d_iF=0, \cD_iy^n=0, \cD_iy_k=0, \cD_ig_{nk}=0\}$,
we respectively obtain  the identities
\be
y_n M^n{}_{ij}=0,
\quad
y^kE_k{}^n{}_{ij}  = -   M^n{}_{ij},   \quad
y_nE_k{}^n{}_{ij}  = g_{kn}M^n{}_{ij}, \quad
E_{mnij}+E_{nmij}=2C_{mnh}
M^h{}_{ij} .
\ee

\ses

In case
of the coefficients
$\{N^k_n,D^k{}_{nm}\}$
 proposed by (2.6) and (2.18)
the direct calculation
of the right-hand parts in (3.2) and (3.3)
results in

 \ses

 \ses

{\bf  Theorem 3.1.} {\it
 The  tensors $M^n{}_{ij}$
 and  $E_k{}^n{}_{ij}$ are represented by the following simple
 and explicit formulas:
\be
M^n{}_{ij}=Fm^nM_{ij}
\ee
and
\be
E_k{}^n{}_{ij}= \Bigl(-l_k m^n+l^n m_k
+
Im_km^n
\Bigr)
M_{ij},
\ee
where
\be
M_{ij}=\D{k_j}{x^i}-\D{k_i}{x^j}.
\ee
}

\ses

\ses

{%\pgbrk}

It proves pertinent to replace in the commutator (3.1)
the partial derivative
$\partial  w^n{}_k/\partial y^h$
by the definition
\be
\cS_h{w^n{}_k}=\D{w^n{}_k}{y^h}+C^n{}_{hk}w^h{}_k-C^m{}_{hk}w^n{}_m
\ee
which has the meaning of the covariant derivative in the tangent
Riemannian space
$\cR_{\{x\}}$.
With
the {\it curvature tensor}
\be
\Rho_k{}^n{}_{ij}=E_k{}^n{}_{ij}
-
M^h{}_{ij}C^n{}_{hk},
\ee
the commutator takes on the form
\be
\lf(\cD_i\cD_j-\cD_j\cD_i\rg) w^n{}_k=
M^h{}_{ij}\cS_h w^n{}_k  -\Rho_k{}^h{}_{ij}w^n{}_h
+\Rho_h{}^n{}_{ij}w^h{}_k.
\ee
The skew-symmetry
\be
\Rho_{mnij}=-\Rho_{nmij}
\ee
holds (cf. the last item in (3.6)).

If we take into account the form of the tensor $C_{ijk}$
indicated in  (2.2), from (3.8) and (3.11) we may conclude that the curvature
tensor is of the following astonishingly simple form:
\be
\Rho_k{}^n{}_{ij}=
(l^n m_k-l_k m^n)
M_{ij}
\equiv
\epsilon^n{}_k
M_{ij}.
\ee
The tensor
$\Rho_{knij}=g_{nm}\Rho_k{}^m{}_{ij}
$
can be represented in the form
\be\Rho_{knij}
=\epsilon_{nk}  M_{ij}.
\ee

We have
$
l^km_n\Rho_k{}^n{}_{ij}=
-
M_{ij}.
$

{%\pgbrk}

\ses

\ses

\setcounter{sctn}{4} \setcounter{equation}{0}

\ses

\ses

\bc
{\bf 4. Riemannian counterparts}
\ec

\ses

\ses

If the Finsler space $\cF_2$ is a Riemannian space with
a Riemannian metric function
$S=\sqrt{a_{ij}y^iy^j}$, where $a_{ij}=a_{ij}(x)$
is a Riemannian metric tensor, we can consider the Riemannian
precursor coefficients
\be
L^k_n={N^k_n}_{\bigl|\text{Riemannian limit}\bigr.}.
\ee
From (2.6)
it follows that
\be
L^k_n
=
-y^k  \fr1S  \D{S}{x^n}
-
Sm^k
\lf(\D{\theta}{x^n}- k_n\rg).
 \ee
On the other hand,
denoting by
$a^k{}_{nh}$
 the Riemannian Christoffel symbols constructed from the
 Riemannian metric  tensor
$a_{mn}$,
we can obtain the equality
\be
a^k{}_{nh}y^h
=
\fr1S\D S{x^n}y^k
+
\lf(\D {\theta}{x^n}+n^h\nabla_n  {\wt b}_h\rg)
Sm^k
\ee
(see (A.9) in Appendix A),
where
$\wt b_h=\wt b_h(x)$ is a vector field chosen to fulfill
\be
\theta(x,\wt b(x))=0
\ee
and the pair $n_i,\wt b_i$ is orthonormal with respect to the tensor
$a_{ij}$.
 The reciprocal pair is  $\{\wt b^i,n^i\}$
with $\wt b^i=a^{ij}\wt b_j$
and
$n^i=a^{ij}n_j$,
where
$a^{ij}$ is the inverse of $a_{ij}$.
The $\nabla_n $ stands for  the Riemannian covariant derivative
taken with
  $a^k{}_{nh}$.
   We get
\be
L^k_{n}=
Sm^k
T_n
-
a^k{}_{nh}y^h\equiv L^k{}_{nh}y^h
\ee
with
\be
L^k{}_{nh}=
-
a^{kj}\epsilon^{\text{Riem}}_{jh}
T_n
-
a^k{}_{nh}
\ee
and
 \be
 T_n=n^h\nabla_n  {\wt b}_h+k_n.
\ee
$\epsilon^{\text{Riem}}_{jh}
=\sqrt{\det(a_{mn})}\,\ga_{jh}$, where
$\ga_{11}=\ga_{22}=0$ and $\ga_{12}=-\ga_{21}=1.$
The metricity property
\be
\D{a_{mn}}{x^i}
+ L^s{}_{im}a_{sn} + L^s{}_{in}a_{ms}
 =0
\ee
holds independently of presence of the vector $T_n$.
In contrast to
the Christoffel symbols $a^k{}_{nh}$,
the coefficients  $L^k{}_{nh}$ obtained are
 not symmetric with respect to the
subscripts.

{%\pgbrk}

Let us take the coefficients
$\bar L^n{}_{ik}=-L^n{}_{ik}$ from (4.6)
and
  construct
  the  tensor
 \be
\bar L_k{}^n{}_{ij} =
\D{\bar L^n{}_{jk}}{x^i}
-
\D{\bar L^n{}_{ik}}{x^j}
+\bar L^m{}_{jk}\bar L^n{}_{im}-\bar L^m{}_{ik}\bar L^n{}_{jm}
\equiv
\bar L_k{}^n{}_{ij}(x).
\ee
We obtain
\be
\bar L_k{}^n{}_{ij} =
\lf(\nabla_i T_j-\nabla_j T_i\rg)
a^{nt}\epsilon^{\text{Riem}}_{tk}
+
a_k{}^n{}_{ij},
\ee
\ses
where
 \be
a_k{}^n{}_{ij} =
\D{a^n{}_{jk}}{x^i}
-
\D{a^n{}_{ik}}{x^j}
+a^m{}_{jk}a^n{}_{im}-a^m{}_{ik}a^n{}_{jm}
\ee
is
the
 Riemannian curvature  tensor
 constructed from the
 Riemannian metric  tensor
$a_{mn}$.
We have taken into account the vanishing
$
\nabla_i \epsilon^{\text{Riem}}_{tk}=0.
$

From the equalities
\be
\nabla_i\wt b^k
=-
n^k
\wt b_m\nabla_in^m, \qquad
\nabla_i n^k
=
-
\wt b^k
n_m\nabla_i\wt b^m
\ee
it follows that
\be
\nabla_i (n^t\nabla_j\wt b_t)
-
\nabla_j (n^t\nabla_i\wt b_t)
=
n^t(\nabla_i \nabla_j\wt b_t-\nabla_j \nabla_i\wt b_t)
=
-
n^t\wt b_la_t{}^l{}_{ij}.
\ee
Therefore,
taking the  $T_i$ from (4.7),
we find that
$$
\nabla_i T_j-\nabla_j T_i=M_{ij}
-
n^t\wt b_la_t{}^l{}_{ij}.
$$
Noting the equality
\be
n^t\wt b_la_t{}^l{}_{ij}
a^{nt}\epsilon^{\text{Riem}}_{tk}=
a_k{}^n{}_{ij}
\ee
(see Appendix A),
we conclude that the tensor
(4.10) can be read merely
\be
\bar L_k{}^n{}_{ij} =
a^{nt}\epsilon^{\text{Riem}}_{tk}
M_{ij}.
\ee
\ses

 If we want to have
$
\bar L^s{}_{ij}=a^s{}_{ij},
$
we must make the choice
\be
k_n=-n^h\nabla_n  {\wt b}_h,
\ee
which entails
 $T_i=0$,  in which case
 the tensor
$\bar L_k{}^n{}_{ij}$ given by  (4.9)
 is the ordinary Riemannian curvature tensor
 $a_k{}^n{}_{ij}$.

{%\pgbrk}

If the Finsler space $\cF_2$ is not a Riemannian space,
it is possible to introduce the
{\it  $\theta$-associated Riemannian space}
 $\cR_{\{\theta\}}$
as follows.

The angle function
$\theta=\theta(x,y)$ is defined
(from  equation (2.3))
up to an arbitrary integration constant
 which may depend on $x$,
which is the reason why  $ d_n\theta$
 should not be put to be zero in (2.8)
(in distinction from the vanishing   $ d_nF=0$).
There exists the freedom to make the redefinition
$\theta\to\theta+C(x)$.
To specify the value of $\theta$ unambiguously  in a fixed
tangent space $T_xM$, we need  in this $T_xM$
an axis
 from which the value
is to be measured. Let the distribution of these axes
over the base manifold
be assigned by means of a
contravariant vector field ${ b}^i={ b}^i(x)$.
 Then we obtain precisely the equality
$
\theta(x,  b(x))=0
$
which does not permit  the  redefinitions  anymore.

\ses

 It is appropriate to construct the
{\it osculating } Riemannian metric tensor
$
a_{mn}(x)=g_{mn}\lf(x, b(x)\rg)
$
and introduce
the normalized vector $\wt b^i=b^i/\sqrt{a_{mn}b^mb^n}$.
Because of the homogeneity,
$g_{mn}\lf(x, b(x)\rg)=g_{mn}\bigl(x,\wt b(x)\bigr)$.
The vector  $n_i(x)$ can be taken to  equal  the value of the derivative
$
\partial {\theta}/\partial {y^i}$
at the argument pair  $\bigl(x,\wt b(x)\bigr)$.
Then,
because of  $ \theta(x,  \wt b(x))=0 $
and
$g_{ij}=l_il_j+m_im_j$
(see (2.2) and (2.3)),
 the pair $\{\wt b_i,n_i\}$ thus introduced
is orthonormal
with respect to the tensor $a_{ij}$ produced by osculation.
This tensor
$a_{ij}$
introduce the Riemannian space
 $\cR_{\{\theta\}}$
 on the base manifold $M$.
We obtain the equalities \be a_{mn} \wt b^m \wt b^n=1, \qquad
a_{mn} n^m n^n=1,
 \qquad   a_{mn} \wt b^m n^n=0,
\ee
and $F\bigl(x,\wt b(x)\bigr)=1$ together with
\be
a_{mn}(x)=g_{mn}\bigl(x,\wt b(x)\bigr),
\quad
\theta(x, \wt b(x))=0,
\quad
\D{\theta}{y^i}\bigl(x,\wt b(x)\bigr)=n_i(x),
\quad
\D{\theta}{y^i}\bigl(x,n(x)\bigr)=-\wt b_i(x).
\ee
The arisen expansion
$
y^m=\wt b\wt b^m+nn^m
$
is convenient to use in many fragments of evaluations.
The last equality in the list (4.18) is explicated from (2.1).

Now we can download in the space $\cR_{\{\theta\}}$ all the
relations (4.2)-(4.16) formulated above in the Riemannian
precursor space. On doing so, we can conclude after comparing
 (4.15) with (3.15)  that
the tensor $\Rho_{knij}=g_{nm}\Rho_k{}^m{}_{ij}$
is factorable, namely
\be
\Rho_{knij}=f_1
\bar L_{knij}  \quad \text{with} ~ ~
\bar L_{knij}= a_{nh}\bar L_k{}^h{}_{ij}
\equiv
\bar L_{knij}(x),
\ee
\ses
where
\be
f_1=
 \sqrt{\fr{\det(g_{hl})}{\det(a_{mn})}}.
\ee

\ses

We have arrived at the following theorem.

\ses

\ses

{\bf  Theorem 4.1.} {\it The curvature tensor
$\Rho_k{}^n{}_{ij}$
of the Finsler space $\cF_2$ equipped with
the
angle-preserving  connection is such that the tensor
$\Rho_{knij}=g_{nm}\Rho_k{}^m{}_{ij}=\Rho_{knij}(x,y)$
is proportional to the tensor
$\bar L_{knij}= a_{nh}\bar L_k{}^h{}_{ij}=\bar L_{knij}(x)$
which does not involve any dependence on tangent vectors.
The factor of proportionality $f_1$ is expressed through the determinants
of metric tensors, according to }
(4.20).

\ses

\ses

{%\pgbrk}

\ses

\ses

\setcounter{equation}{0}

\bc
{ \bf 5. Conclusions}
\ec

\ses

\ses

In the Riemannian geometry
the contraction $a^k{}_{nh}y^h$ of the
Christoffel symbols
 $a^k{}_{nh}$
with the tangent vector $y$
admits the angle representation (4.3)--(4.4):
\be
a^k{}_{nh}y^h
=
\fr1Sy^k
\D S{x^n}
+
\lf(\D {\theta}{x^n}-t_n(x)\rg)
Sm^k,
\ee
where
$ t_n(x)=-n^h\nabla_n  {\wt b}_h$
and
$\theta(x,\wt b(x))=0$.
Why don't lift  the representation to the Finsler level
to take the coefficients $N^k_n=N^k_n(x,y)$
in the operator
$
d_n=\partial_{x^n}+N^k_n(x,y)\partial_{y^k}
$
to be of the similar form?
Our proposal in (2.6) was of this kind,
namely,
\be
N^k_n
=
-\fr1F y^k    \D{F}{x^n}
-
\lf(\D{\theta}{x^n}- k_n(x)\rg)F  m^k.
 \ee
\ses
At any $k_n$, the vanishing
$d_nF=d_n\theta-k_n=d_n(\theta_2-\theta_1)=0$
immediately ensues from this proposal.
It is a big (and good) surprise that the vanishing
$y_kN^k{}_{nmi}=0$ ensues also, which enables us
to obtain the covariant derivative
$\cD_n$ possessing the metric property $\cD_ng_{ij}=0$,
where
$\cD_ng_{ij}= d_ng_{ij} - D^h{}_{ni}g_{hj}- D^h{}_{nj}g_{ih}$
with the connection coefficients
$D^h{}_{ni}=-N^h{}_{ni}$. If we want to obtain torsionless coefficients
in the Riemannian limit of these $D^h{}_{ni}$, we should take the vector
field $k_n(x)$ to be the  field
$ t_n(x)=-n^h\nabla_n  {\wt b}_h$
 which enters the Riemannian version (5.1).

The induced
parallel transports of the objects
$\{F,\,\theta_2-\theta_1, \,g_{ij}\}$
along the horizontal curves (running on  the base manifold $M$)
are
represented
infinitesimally  by
the elements
$\{dx^nd_nF,\,dx^nd_n(\theta_2-\theta_1),\,
dx^n\cD_ng_{ij}\}$
which all are  the naught
because of $d_nF=d_n(\theta_2-\theta_1)=\cD_ng_{ij}=0$.
Therefore
the transports
 realize  isometries of the tangent Riemannian spaces
 $\cR_{\{x\}}$ supported by points $x\in M$,
 taking indicatrices into indicatrices.
The coefficients  $ N^k_n$ given by (2.6)
 are in general non-linear with respect to the variable $y$.

In the Riemannian case, the right-hand part of (5.1)
 can be expressed through the
 Christoffel symbols and, therefore, can be constructed from the first
 derivatives of the metric tensor.
 This is the privilege of the Riemannian geometry which lives
 in the ground floor of the Finsler building,
 ---
the right-hand part of the Finsler coefficients (5.2)
is not a composition of  partial derivatives of the Finsler metric tensor.
In distinction to the Riemannian geometry which provides us with simple
and explicit angle (see (A.8) in Appendix A),
the Finsler angle function $\theta=\theta(x,y)$ is defined by the
partial differentiable equation
 (2.3)
which cannot be integrated explicitly, except for rare
particular  cases
of the Finsler metric function.

\ses

The second big surprise is that the angle-preserving connection
obtained in this way
admits
the $C^{\infty}$-regular realization  globally
regarding
the dependence on tangent vectors.
Such a realization takes place for the Finsleroid-regular metric
function, the Randers metric function,
and probably for many other Finsler  metric functions.

The Finsler connection obtained does not need any facility which could be provided
by the geodesic spray coefficients. Due attention to the angle wisdom
is sufficient: the Finsler space is connected by its angle structure,
 similarly to the well known property of Riemannian geometry.

Our consideration was restricted by the dimension $2$.
Development of due
extensions to higher dimensions is the problem
of urgent kind.

\ses

{%\pgbrk}

\setcounter{equation}{0}

\bc
{ \bf Appendix A: Involved evaluations}
\ec

\ses

\ses

Let us verify  theorem 3.1.
Using
$d_im^k=N^k{}_{ih}m^h$
(see (2.22)) and (2.11), we obtain
$$
d_im^k=
-
m^hl^k \D{l_h}{x^i}
+
(Im^k + l^k)\breve P_i
-
m^h m^k
\D{m_h}{x^i}.
$$
\ses
Also,
$$
d_i\D{\theta}{x^j}-d_j\D{\theta}{x^i}=
N_i^t\D{\fr1Fm_t}{x^j}    -   N_j^t\D{\fr1Fm_t}{x^i}
$$

\ses

$$
=
N_i^t
\lf[
 -\fr1{F^2} \D F{x^j}m_t+\fr1F\D{m_t}{x^j}
\rg]
-
N_j^t
\lf[
 -\fr1{F^2}\D F{x^i}m_t+\fr1F\D{m_t}{x^i}
\rg]
$$

\ses

   \ses

$$
=
  m^t\breve P_i           \lf[   \fr1{F}\D F{x^j}m_t-\D{m_t}{x^j}  \rg]
 -    m^t\breve P_j          \lf[   \fr1{F}\D F{x^i}m_t-\D{m_t}{x^i}  \rg].
 $$
\ses
Using these  formulas in (3.2) and
keeping also in mind that
$ d_iF=0$ and $  d_iy^k=N^k_i ,$
from (2.6) we get
$$
M^n{}_{ij}=
-N^n_i \fr1F\D F{x^j}+N^n_j \fr1F\D F{x^i}
-l^n
N^h_i
          \D{\D F{x^j}}{y^h}
+l^n  N^h_j
           \D{\D F{x^i}}{y^h}
$$

\ses

$$
+
m^hy^n
\lf(\breve P_j\D{l_h}{x^i}-\breve P_i\D{l_h}{x^j}\rg)
+
Fm^h m^n
\lf(\breve P_j\D{m_h}{x^i}-\breve P_i\D{m_h}{x^j}\rg)
$$

\ses

$$
-Fm^n
  \lf(
  m^t\breve P_i           \lf(   \fr1{F}\D F{x^j}m_t-\D{m_t}{x^j}  \rg)
 -    m^t\breve P_j           \lf(   \fr1{F}\D F{x^i}m_t-\D{m_t}{x^i}  \rg)
   \rg  )
+Fm^nM_{ij},
$$
\ses
which reduces to
$$
M^n{}_{ij}=
-l^n
N^h_i
          \D{\D F{x^j}}{y^h}
+l^n  N^h_j
           \D{\D F{x^i}}{y^h}
+
m^h
y^n
\lf(\breve P_j\D{l_h}{x^i}-\breve P_i\D{l_h}{x^j}\rg)
+Fm^nM_{ij}
=Fm^nM_{ij}.
$$
The representation (3.7) holds.
Differentiating  (3.7) with respect of $y^k$ and taking into account
the formula (3.5)
together with (2.12), the representation (3.8) is obtained.
The theorem is valid.

{%\pgbrk}

\ses

Let us verify the equality (4.14).
Since the dimension is $N=2$, the Riemannian curvature tensor
$a_{tlin}=a_{lh}a_{t}{}^h{}_{in}$ possesses the representation
\be
a_{tlij}=-R(a_{tj}a_{li}-  a_{ti}a_{lj} ).
\ee
The meaning of the scalar $R$ is explained by the equality
\be
a^{ti}a^{lj}a_{tlij}=2R.
\ee
Inserting the expansion
$
a_{tn}=\wt b_t\wt b_n+n_tn_n
$
in (A.1) results in the factorization
\be
a_{tlij}=-R(n_t\wt b_l-n_l\wt b_t)(n_j\wt b_i-n_i\wt b_j).
\ee
It follows that
$ n^t\wt b^la_{tlij}=-R(n_j\wt b_i-n_i\wt b_j)$
and therefore
\be
 n^t\wt b^la_{tlij}({\wt b}_kn_n-n_k{\wt b}_n)
 =
 -R(n_j\wt b_i-n_i\wt b_j)({\wt b}_kn_n-n_k{\wt b}_n)
 =-a_{knij}.
\ee
We can substitute here
$\epsilon^{\text{Riem}}_{kn}$ with
$
{\wt b}_kn_n-n_k{\wt b}_n.$
The equality (4.14) is valid.

{%\pgbrk}

\ses

In the remainder we verify
the representation of $a^k{}_{nh}y^h$ indicated in (2.25).

In the Riemannian limit we have $||b||=1$ and may use the notation
$b_i$ instead of $\wt b_i$.
We start with the  Riemannian representations
$$
 a_{ij}=n_in_j+b_ib_j, \qquad l^i=\fr{y^i}S, \qquad l_i=a_{ij}l^j, \qquad
S=\sqrt{a_{mn}y^my^n},
$$
\ses$$
l^i=n^i\sin\theta+b^i\cos\theta,  \qquad
l_i=n_i\sin\theta+b_i\cos\theta,
$$
and
$$
m^i=n^i\cos\theta-b^i\sin\theta,
  \qquad
m_i=n_i\cos\theta-b_i\sin\theta.
$$
\ses
We have
$$
n=S\sin\theta, \qquad b=S\cos\theta,
$$
which entails
\be
m^i=\fr1S(bn^i-nb^i), \quad m_i=\fr1S(bn_i-nb_i),
\quad l^il_i=m^im_i=1, \quad l^im_i=0.
\ee

\ses

We find
$$
 a_{ij}-l_il_j
 =n_in_j+b_ib_j
 -
n_i\sin\theta
(n_j\sin\theta+b_j\cos\theta)
-
b_i\cos\theta
(n_j\sin\theta+b_j\cos\theta)
$$

\ses

$$
 =n_in_j\cos^2\theta+b_ib_j\sin^2\theta -(n_ib_j+n_jb_i)\sin\theta\cos\theta,
$$
\ses
or
$$
 a_{ij}-l_il_j
 =(n_i\cos\theta-b_i\sin\theta)
(n_j\cos\theta-b_j\sin\theta)\equiv m_im_j.
$$

Noting
that
$$
\D{l_i}{x^n}=
\fr1S
(nn_{i,n}+b b_{i,n} +m_i \,\theta_{,n}) ,\quad
\D{m_i}{x^n}=
 \fr1S(bn_{i,n}-n b_{i,n} - y_i\theta_{,n})
 $$
(by $\{,n\}$ we denote the partial derivative with respect to $x^n$),
from the expansion
$
 a_{ij}=l_il_j+m_im_j
$
we obtain
$$
 a_{ij,n}=
\fr1S
(nn_{i,n}+b b_{i,n} +Sm_i \,\theta_{,n})l_j
+
\fr1S
(nn_{j,n}+b b_{j,n} +Sm_j \,\theta_{,n})
l_i
$$

\ses

$$
+
 \fr1S(bn_{i,n}-n b_{i,n} - y_i\theta_{,n})  m_j
+
 \fr1S(bn_{j,n}-n b_{j,n} - y_j\theta_{,n})  m_i,
 $$
\ses
which can be written in the simpler form
\be
 a_{ij,n}=
\fr1S
(nn_{i,n}+b b_{i,n})l_j
+
\fr1S
(nn_{j,n}+b b_{j,n})
l_i
+
 \fr1S(bn_{i,n}-n b_{i,n})  m_j
+
 \fr1S(bn_{j,n}-n b_{j,n})  m_i.
\ee

{%\pgbrk}

Let us compose the sum
$$
 a_{in,j}+ a_{jn,i}  -  a_{ij,n}=
$$

\ses

$$
\fr1S
(nn_{i,j}+b b_{i,j})l_n
+
\fr1S
(nn_{n,j}+b b_{n,j})
l_i
+
 \fr1S(bn_{i,j}-n b_{i,j})  m_n
+
 \fr1S(bn_{n,j}-n b_{n,j})  m_i
 $$

\ses

\ses

$$
+
\fr1S
(nn_{j,i}+b b_{j,i})l_n
+
\fr1S
(nn_{n,i}+b b_{n,i})
l_j
+
 \fr1S(bn_{j,i}-n b_{j,i})  m_n
+
 \fr1S(bn_{n,i}-n b_{n,i})  m_j
 $$

\ses

\ses

$$
-
\fr1S
(nn_{i,n}+b b_{i,n})l_j
-
\fr1S
(nn_{j,n}+b b_{j,n})
l_i
-
 \fr1S(bn_{i,n}-n b_{i,n})  m_j
-
 \fr1S(bn_{j,n}-n b_{j,n})  m_i,
 $$
\ses
from which we obtain
$$
( a_{in,j}+ a_{jn,i}  -  a_{ij,n})y^j=
$$

\ses

$$
=
\fr1S
(nn_{i,j}+b b_{i,j})  y^jl_n
+
\fr1S
(nn_{n,j}+b b_{n,j})
y^jl_i
+
 \fr1S(bn_{i,j}-n b_{i,j}) y^j m_n
+
 \fr1S(bn_{n,j}-n b_{n,j}) y^j m_i
 $$

\ses

\ses

$$
+
\fr1S
(nn_{j,i}+b b_{j,i})y^jl_n
+
(nn_{n,i}+b b_{n,i})
+
 \fr1S(bn_{j,i}-n b_{j,i}) y^j m_n
 $$

\ses

\ses

$$
-
(nn_{i,n}+b b_{i,n})
-
\fr1S
(nn_{j,n}+b b_{j,n})
y^jl_i
-
 \fr1S(bn_{j,n}-n b_{j,n}) y^j m_i,
 $$
\ses
or
$$
( a_{in,j}+ a_{jn,i}  -  a_{ij,n})y^j=
\fr1S
(nn_{i,j}+b b_{i,j})  y^jl_n
+
\fr1S
\Bigl(n(n_{n,j}-n_{j,n})  +b (b_{n,j}-b_{j,n})
\Bigr)
y^jl_i
$$

\ses

\ses

$$
+
 \fr1S(bn_{i,j}-n b_{i,j}) y^j m_n
+
 \fr1S
 \Bigl(b
 (n_{n,j}-n_{j,n})    -n (b_{n,j}-b_{j,n})\Bigr) y^j m_i
 $$

\ses

\ses

\be
+
\fr1S
(nn_{j,i}+b b_{j,i})y^jl_n
+
n(n_{n,i}-n_{i,n})  +b (b_{n,i}-b_{i,n})
+
 \fr1S(bn_{j,i}-n b_{j,i}) y^j m_n.
\ee
\ses
We use here the equalities
$$
\nabla_nb_k=p_nn_k, \qquad \nabla_nn_k=-p_nb_k, \qquad
\text {with}
 ~ ~ p_n=n^h\nabla_nb_h,
$$
where
 $\nabla_n $ stands for  the Riemannian covariant derivative
 taken with
 the Riemannian Christoffel symbols
 $a^k{}_{nh}$
 constructed from the
 Riemannian metric  tensor
$a_{mn}$,
which entails
$$
n(n_{n,j}-n_{j,n})  +b (b_{n,j}-b_{j,n})=S(p_jm_n-p_nm_j).
$$

{%\pgbrk}

\nin
The method yields
$$
( a_{in,j}+ a_{jn,i}  -  a_{ij,n})y^j=
\fr1S
(nn_{i,j}+b b_{i,j})  y^jl_n
+
y^jp_jm_n
l_i
$$

\ses

\ses

$$
+
 \fr1S(bn_{i,j}-n b_{i,j}) y^j m_n
+
\Bigl[b(-p_jb_n+p_nb_j)-n(p_jn_n-p_nn_j)
\Bigr]
l^j  m_i
 $$

\ses

\ses

$$
+
\fr1S
(nn_{j,i}+b b_{j,i})y^jl_n
+
S(p_im_n-p_nm_i)
+
 \fr1S(bn_{,i}-n b_{,i})  m_n
 $$

\ses

\ses

\ses

$$
=\fr1S
\Bigl(n(n_{i,j}+n_{j,i})+b (b_{i,j}+b_{j,i})\Bigr)  y^jl_n
+
y^jp_jm_n
l_i
$$

\ses

\ses

$$
+
 \fr1S
 \Bigl(b(-p_jb_i+p_ib_j)
 -n (p_jn_i-p_in_j)
 \Bigr) y^j m_n
+
Sp_im_n
+
 \fr2S(bn_{,i}-n b_{,i})  m_n
-y^jp_jl_nm_i
 $$

\ses

             \ses

\ses

$$
\!=\!\fr1S
\Bigl(2SS_{,i}+n(-p_jb_i+p_ib_j)  y^j +b (p_jn_i-p_in_j)  y^j\Bigr)l_n
\!+\!
2Sp_im_n
\!+\!
 \fr2S(bn_{,i}\!-\!n b_{,i})  m_n
\!-y^jp_jl_nm_i.
 $$
\ses
We obtain eventually the expansion
$$
\fr12
( a_{in,j}+ a_{jn,i}  -  a_{ij,n})y^j=
S_{,i}l_n
+
\Bigl(p_i+
 \fr1{S^2}(bn_{,i}-n b_{,i})
 \Bigr) S  m_n.
 $$
\ses
Since
\be
\theta=\arctan\fr nb, \qquad
\theta_{,n}=\fr1{S^2}(bn_{,n}-nb_{,n}),
\ee
\ses
we get
$$
\fr12
( a_{in,j}+ a_{jn,i}  -  a_{ij,n})y^j=
S_{,i}l_n
+
(\theta_{,i}+p_i)
 S  m_n.
$$
\ses
 Thus we have
\be
a^n{}_{ij}y^j
=
S_{,i}l^n
+
(\theta_{,i}+p_i)
 S  m^n
\ee
which is the representation (2.25).

From (A.8) we may observe that
$$
\theta(x,b(x))=0.
$$

{%\pgbrk}

\ses

\ses

\setcounter{equation}{0}

\bc
{ \bf Appendix B: Finsleroid-regular space
with connection}
\ec

\ses

\ses

We start with the notion of a 2-dimensional Riemannian space
$\cR_2=(M,{\cal S})$ to specify an attractive
 2-dimensional
Finsler space  over  $\cR_2$.
We shall assume that
the  background two-dimensional manifold $M$ admits
introducing two linearly independent
covariant vector fields, to be presented  by
the  non--vanishing 1-forms $\wt b=\wt b(x,y)$ and  $n=n(x,y)$.
With respect to  natural local coordinates $\{x^n\}$ in the manifold $M$
the expansions
\be
\wt b=\wt b_i(x)y^i, \qquad n=n_i(x)y^i, \qquad
 S= \sqrt{a_{ij}(x)y^iy^j},
\ee
represent  the 1-forms and the Riemannian metric $\cS$, with
$a_{ij}$ standing for the covariant components of the Riemannian metric tensor of the
space  $\cR_2$.  The contravariant components  $a^{ij}$ are defined by means
of the reciprocity  $a^{in}a_{nj}=\de^i{}_j$.
The covariant index of the vectors $\wt b_i$  and
$ n_i$
will be raised by means of the Riemannian rule
$ \wt b^i=a^{ij}\wt b_j,\,  n^i=a^{ij}n_j,$
 which inverse reads $ \wt b_i=a_{ij}\wt b^j, \, n_i=a_{ij}n^j.$
Also, we assume that the vectors introduced  are  orthonormal
relative to the Riemannian metric ${\cal S}$, that is,
\be
a_{mn} \wt b^m \wt b^n=1, \qquad   a_{mn} n^m n^n=1,
 \qquad   a_{mn} \wt b^m n^n=0.
\ee

{%\pgbrk}

Let a positive scalar $c=c(x)$ be given which is ranged as follows:
\be
 0 < c < 1.
 \ee
We are entitled to construct
 the 1-form
 \be
 b=b_i(x)y^i=c\wt b,
 \ee
so that
\be
\wt b_i=\fr1cb_i, \qquad  \wt b^i=\fr1cb^i, \qquad
c=|| b||_{\text{Riemannian}} \equiv \sqrt{ a^{mn}  b_m b_n},
\ee
where $  b^i=a^{ij} b_j$.

Since $c<1$,  we get
\be
S^2-b^2>0 ~~ ~~ \text{whenever} ~ ~ y\ne 0
\ee
and may conveniently use the variable
\be
q:=\sqrt{S^2-b^2}
\ee
which {\it does not vanish anywhere on}
 $ T_xM\backslash 0$.
 Obviously, the inequality
\be
q^2 \ge \fr{1-c^2}{c^2}\,b^2
\ee
is valid.

{%\pgbrk}

We have
\be
q=\sqrt{\fr{1-c^2}{c^2}b^2+n^2},
\ee
so that
\be
q\D q{y^m}=\fr{1-c^2}{c^2}bb_m+nn_m
\ee
\ses
and
$$
\D{\fr bq}{y^m}=\fr 1q b_m-\fr b{q^3}\lf(\fr{1-c^2}{c^2}bb_m+nn_m\rg),
$$
or
\ses\\
\be
\D{\fr bq}{y^m}=   - \fr n{q^3}(bn_m-nb_m).
\ee

\ses

{%\pgbrk}

The space described below involves the characteristic parameter
\be
g=g(x)\in(-2,2).
\ee
It is convenient to introduce the quantities
\be
h=\sqrt{1-\fr{g^2}4}, \qquad G=\fr gh.
\ee
\ses
The functions $B=B(x,y)$ and  $B_1=B_1(x,y)$  are specified
as follows:
\be
B =  b^2+gbq+q^2, \qquad B_1=\fr1{c^2}(b+gc^2q), \qquad
B-bB_1=n^2.
\ee
\ses
The function
$B $   is {\it   positively definite}, for
$$
b^2+gqb+q^2
=
\fr 12
\Bigl[(b+g_+q)^2+(b+g_-q)^2\Bigr],
$$
where $ g_+=(1/2)g+h$ and $ g_-=(1/2)g-h$.
In the limit $g\to 0$,
the function $B$  degenerates to the
 quadratic form  of the input Riemannian metric tensor:
\be
 B|_{_{g=0}}=b^2+q^2 \equiv S^2.
\ee
Also,
\be
n|_{_{y^n=b^n}}=0, \quad
b|_{_{y^n=b^n}}=c^2, \quad
q|_{_{y^n=b^n}}=c\sqrt{1-c^2},
\quad S^2|_{_{y^n=b^n}}=c^2,
\quad \eta B|_{_{y^n=b^n}}=c^2,
\ee
where
\be
\eta=\fr 1{1+gc\sqrt{1-c^2}}.
\ee
It can readily be verified that on the definition range $g\in(-2,2)$ of the parameter
 $g$  we have
$
\eta>0.
$

{%\pgbrk}

\ses

The basic metric function  $K$ comes  from  the following definition.

\ses

\ses

 {\it Key  Definition}.
  The scalar function $K(x,y)$ given by the formulas
\be
K(x,y)=
\sqrt{B(x,y)}\,J(x,y),  \qquad
J(x,y)=\e^{-\frac12G(x)f(x,y)},
\ee
where
\be
f=
-\arctan \fr G2+\arctan\fr{L}{hb},
\qquad  {\rm if}  \quad b\ge 0,
\ee
and
\be
f= \pi-\arctan\fr G2+\arctan\fr{L}{hb},
\qquad  {\rm if}
 \quad b\le 0,
\ee
 with
 \be
 L = q+\fr g2b,
\ee
\ses\\
is called
the {\it  Finsleroid-regular  metric function}.

\ses

\ses

This metric function $K$ is not absolute homogeneous.

\ses

The function $L$  obeys   the identity
\be
L^2+h^2b^2=B.
\ee

\ses

\ses

\ses

 {\large  Definition}.  The arisen  Finsler space
\be
\cF\cR^{PD}_{g;c} :=\{\cR_2;\,b_i(x);\,g(x);\,K(x,y)\}
\ee
is called the
 {\it Finsleroid-regular  space}.

\ses

\ses

The upperscripts $PD$ mean that the space is positive-definite.

\ses

\ses

 {\large  Definition}.  Within  any tangent space $T_xM$, the
  metric function $K(x,y)$
 produces the $\cF\cR^{PD}_{g;c}$-{\it circle}
 \be
 \cF\cR^{PD}_{g;c\, \{x\}}:
 =\{y\in   \cF\cR^{PD}_{g;c\, \{x\}}: y\in T_xM , K(x,y)\le 1\}.
  \ee

\ses

 \ses

 {\large  Definition}. The
 $\cF\cR^{PD}_{g;c}$-{\it indicatrix}
 $ I\cR^{PD}_{g;c\, \{x\}} \subset T_xM$
 is the boundary of the $\cF\cR^{PD}_{g;c}$- circle,
  that is,
 \be
 I\cR^{PD}_{g;c\, \{x\}} :=\{y\in I\cR^{PD}_{g;c\, \{x\}} : y\in T_xM, K(x,y)=1\}.
  \ee

\ses

 \ses

 {\large  Definition}.
 The scalar $g(x)$ is called the {\it Finsleroid charge}.
The 1-form $b=b_i(x)y^i$ is called
 the $\cF\cR^{PD}_{g;c}$-{\it axis one-form}.

\ses

\ses

The  indicatrix (B.25) is obviously symmetric under reflections with respect to
the direction assigned by the vector $b_i$.

The metric function $K$ given by (B.18)-(B.21)
 is regular of the class $C^{\infty}$ regarding the $y$-dependence.
Formally,
this high regularity comes from the circumstance that $q>0$ and the derivatives
$\partial q/\partial y^m$
and
$\partial(b/q)/\partial y^m$ (see (B.10) and (B.11)) are of this class $C^{\infty}$.
The  made assumption $c\in (0,1)$ is essential, for at $c=1$ the quantity $q$ may vanish
(at $n=0$, that is, on the axis of the indicatrix).

{%\pgbrk}

We shall meet the function
\be
\nu:=q+(1-c^2)gb
\ee
for which
\be
\nu >0 \quad \text{when} \quad |g|<2.
\ee
Indeed, if $gb>0$, then  the right-hand part of (B.26) is positive.
 When $gb<0$, we may note that
 at any fixed $c$ and $b$ the minimal value of $q$ equals
$\sqrt{1-c^2}|b|/c$ (see (B.8)), arriving again at (B.27).

The identities
\be
\fr{c^2S^2-b^2}{q\nu}=1-(1-c^2)\fr B{q\nu},
\qquad
gb(c^2S^2-b^2)=qB-\nu S^2
\ee
are of great help to take into account when performing evaluations of key objects entailed.

\ses

We straightforwardly obtain the  representations
\be
y_i~:=\fr12\D {K^2}{y^i}=( u_i+gqb_i ) \fr{K^2}B, \qquad u_i=a_{ij}y^j,
\ee
\ses
and
\be
g_{ij}~:=\D{y_i}{y^j}=
\Biggl[a_{ij}
+\fr g{B}\Biggl(\lf(gq^2-\fr{bS^2}q\rg)b_ib_j-\fr bqu_iu_j+
\fr{ S^2}q(b_iu_j+b_ju_i)\Biggr)\Biggr]\fr{K^2}B,
\ee
together with the reciprocal components
\be
g^{ij}=
\Biggl[a^{ij}+\fr g{\nu}(bb^ib^j-b^iy^j-b^jy^i)+\fr g{B\nu}(b+gc^2q)y^iy^j
\Biggr]\fr B{K^2}.
\ee
\ses

The associated Riemannian metric tensor $a_{ij}$ has the meaning
\be
 a_{ij}=g_{ij}\bigl|_{g=0}\bigr. .
\ee

\ses

The determinant of the Finsler metric tensor presented by (B.30)  is everywhere positive:
\be
\det(g_{ij})=\fr {\nu}q\biggl(\fr{K^2}B\biggr)^2\det(a_{ij})>0.
\ee

\ses

Using the function $X$ given by
\be
\fr1X=    2+\fr{(1-c^2)B}{q\nu} \equiv 3-\fr{c^2n^2}{q\nu},
\ee
\ses
we obtain
\be
A_i~:=K
\D{\ln\lf(\sqrt{\det(g_{mn})}\rg)}{y^i}=\fr {Kg}{2qB}
\fr1X
(S^2b_i- bu_i)
\equiv
\fr {Kg}{2qB}
\fr nX
(nb_i- bn_i)
\ee
\ses
and
\be
A^i~:=g^{ij}A_j=
\fr {g}{2XK\nu}
\Bigl[ Bb^i
-   (b+gqc^2)y^i
\Bigr]
\equiv
\fr {g}{2XK\nu}
( Bb^i- c^2B_1y^i),
\ee
\ses
which entails the contraction
\be
A^iA_i=
\fr{g^2}{4X^2}
\lf(3-\fr1X\rg)
\equiv
\fr{g^2}{4X^2}\fr{c^2n^2}{q\nu}.
\ee

{%\pgbrk}

With the function
\be
 T=\fr1{c}
\sqrt{\fr{\nu} q}\fr{K^2}{B}
\equiv
\fr1c
\sqrt{\fr{\det(g_{ij})}{\det(a_{mn})}},
\ee
\ses
we introduce the vector
\be
  m_n= \fr1K  T  (bn_n-nb_n)
  \ee
\ses
which coincides with the vector (2.1).
The contravariant components $m^n=g^{nk}m_k$ are found with the help of (B.31)
to read
$$
  m^n=
c
\sqrt{\fr q{\nu}}
\lf(-
\fr1{c^2}
n
b^n
+
\lf(\fr1{c^2}b+gq\rg)
 n^n
\rg)
\fr1{K},
$$
\ses
or
\be
  m^n=
c
\sqrt{\fr q{\nu}}
\lf(
\fr1{c^2}
(bn^n-nb^n)
+
gq
 n^n
\rg)
\fr1{K}.
\ee
We can write also
$$
bm^n=c \sqrt{\fr q{\nu}}
(Bn^n-ny^n)\fr1K.
$$
The unit norm equality
$$
m^nm_n=1
$$
is valid.

\ses

Comparing (B.35) and (B.37) with (B.39) shows that
$$
\fr{A_k}{\sqrt{A_hA^h}}=
m_k\ga \quad \text{whenever} ~ g\ne0 ~~ \text{and} ~ n\ne0,
$$
with
$$
\ga=
-
\sign(g)\sign(n),
$$
where
the designation $`\sign`$ stands for the function: $\sign(x)=1$, if $x>0$, and
$\sign(x)=-1$, if $x<0$.

From (B.35) and (B.39) the representation
\be
A_i=Im_i
\ee
ensues,
where
\be
I= -gc\fr1{2X}\fr n{\sqrt{q\nu }}.
\ee

We shall use the angle function $\theta$ defined by the equation
\be
K\D{\theta}{y^n}=m_n,
\ee
where the components $m_n$ placed in the right-hand part are taken from (B.39).

The $I$ given by (B.42) and
the components $A_i=A_i(x,y)$
 reveal the vanishing at the values $y^n=\pm b^n(x)$:
$$
I(x,b(x))=I(x,-b(x))=0,  \qquad A_i(x,b(x))=A_i(x,-b(x))=0,
$$
that is, on the indicatrix axis
(and nowhere else), because they involve the factor $n$.

At the same time, on  all the slit tangent bundle $T_xM\backslash 0$
the vector $m_k$ introduced in (B.39)
 is
 $C^{\infty}$-regular  regarding
the dependence on tangent vectors.

       \ses

   {%\pgbrk}

Since
$$
K|_{_{y^n=b^n}}= c\fr1{\sqrt{\eta}}\e^{-\frac12Gj}
$$
with
$$
j=  -\arctan \fr G2+\arctan\fr{\sqrt{1-c^2}+\fr12gc}{hc},
$$
where we have used  (B.18) and (B.19),
\ses
it may be convenient to rescale the metric function $K$ as follows:
\be
K\to\breve K= \sqrt{\eta}\e^{\frac12G j}K.
\ee
We obtain
the  function
\be
\breve K(x,y)=
\sqrt{\eta B(x,y)}\,\breve J(x,y),  \qquad
\breve J(x,y)=\e^{-\frac12G(x)\breve f(x,y)},
\ee
where
\be
\breve f=
\arctan\fr{L}{hb}  -\arctan\fr{\sqrt{1-c^2}+\fr12gc}{hc},
\qquad  {\rm if}  \quad b\ge 0,
\ee
and
\be
\breve f= \pi+\arctan\fr{L}{hb}-\arctan\fr{\sqrt{1-c^2}+\fr12gc}{hc},
\qquad  {\rm if}
 \quad b\le 0,
\ee
which possesses the property
\be
\breve K|_{_{y^n=b^n}}= c.
\ee
In terms of the metric function normalized in this way the norm of the input
vector $b^i$ is the same with
respect to the considered Finsler space
as well as the underlined  Riemannian space $\cR_2$, that is,
\be
\breve K(x,b(x))=\sqrt{a_{ij}(x)b^i(x)b^j(x)}=c(x).
\ee

{%\pgbrk}

From (B.30)  we get
$$
g_{ij}|_{_{y^n=b^n}}=
\biggl[a_{ij}
+ g\eta
\Bigl(gc^2(1-c^2)
-\fr{c^4}{c\sqrt{1-c^2}}
-\fr{c^2}{c\sqrt{1-c^2}}
+  2\fr{c^2}{c\sqrt{1-c^2}}
\Bigr)
\wt b_i\wt b_j
\biggr]
\fr{\eta}{c^2}
K^2|_{_{y^n=b^n}}
$$

\ses

\ses

$$
=
\biggl[a_{ij}
+ g\eta
\Bigl(gc^2(1-c^2)
+c\sqrt{1-c^2}
\Bigr)
\wt b_i\wt b_j
\biggr]
\fr{\eta}{c^2}
K^2|_{_{y^n=b^n}},
$$
\ses
or
$$
g_{ij}|_{_{y^n=b^n}}=
\lf(
a_{ij}
+ gc\sqrt{1-c^2}\,
\wt b_i\wt b_j
\rg)
\fr{\eta}{c^2}
K^2|_{_{y^n=b^n}}.
$$
\ses
Introducing the osculating tensor
\be
\breve a_{ij}~:=g_{ij}|_{_{y^n=b^n}}
\ee
leads therefore to
\be
a_{ij}
=
\e^{Gj}\breve a_{ij}
-
 gc\sqrt{1-c^2}\,
\wt b_i\wt b_j.
\ee
We may make here the expansion
$$
\breve a_{ij}
=\breve b_i\breve b_j+\breve n_i\breve n_j,
$$
obtaining
\be
\breve b_i=\fr1{\sqrt{\eta}}\e^{-\frac12Gj}\wt b_i, \qquad
  \breve   n_i= \e^{-\frac12Gj}n_i.
\ee
\ses
In terms of the  new variables
$
\{\breve n=\breve n_iy^i, ~   \breve b=\breve b_iy^i
\}$
we can replace the initial quantities
$
q=\sqrt{ S^2-b^2}
$
 and
 $
S=\sqrt{ a_{ij}y^iy^j}
 $
by the quantities
$
\breve q=\sqrt{ \breve S^2-\breve b^2}
$
and
$
\breve S=\sqrt{ \breve a_{ij}y^iy^j}
 $
which are adaptable to the treatment of
the osculating tensor $\breve a_{ij}$
to be
 the metric tensor of the background
Riemannian space.
In this vein the derivative values
\be
\D{\breve \theta}{y^i}\lf(x,\breve b(x)\rg)=\breve n_i(x), \qquad
\D{\breve \theta}{y^i}\lf(x,-\breve b(x)\rg)=-\breve n_i(x)
\ee
are obtained
(cf. (4.18)).

{%\pgbrk}

Now, differentiating $K  m^k$
taken from (B.40) yields
\ses
$$
\D{K  m^k}{y^m}=
- Km^k\fr{q}{2\nu}(1-c^2)g\D {\fr bq}{y^m}
+c
\sqrt{\fr q{\nu}}
\lf(
\fr1{c^2}(b_mn^k-n_mb^k)
+
g\D{q}{y^m}
 n^k
 \rg),
$$
\ses
so that we can write
$$
\D{K  m^k}{y^m}-l_mm^k+l^km_m=
 Km^k\fr{q}{2\nu}(1-c^2)g    \fr n{q^3}\Upsilon_m
$$

\ses

\ses

$$
+c
\sqrt{\fr q{\nu}}
\lf(
\fr1{c^2}\Upsilon^k_m
+
g\fr1q         \lf(\fr{1-c^2}{c^2}bb_m+nn_m\rg)
 n^k
  \rg)
-l_mm^k+l^km_m.
$$
\ses
We use the notation
$$
\Upsilon_m=
bn_m-nb_m, \qquad
\Upsilon^k=
bn^k-nb^k, \qquad
\Upsilon^k_m=
b_mn^k-n_mb^k.
$$
\ses
We straightforwardly have
\ses

$$
-l_mm^k+l^km_m= \fr1B  c
\sqrt{\fr q{\nu}}\, p^k_m
$$
\ses
with
$$
p^k_m=\fr {\nu}q\,
y^k \fr1{c^2}
\Upsilon_m
 -
\lf(
\fr1{c^2}
\Upsilon^k
+
gq
 n^k
\rg)
(u_m+gqb_m)
$$

\ses

\ses

\ses

$$
=
  bb^k \fr1{c^4}     (bn_m-nb_m)   +    nn^k \fr1{c^2}     (bn_m-nb_m)
+\fr{1-c^2}{c^4}
g\fr bq bb^k    \Upsilon_m
+(1-c^2)g\fr bq nn^k \fr1{c^2} \Upsilon_m
$$

\ses

$$
 - \lf(\fr1{c^2}(bn^k-nb^k)+gq n^k\rg)\fr b{c^2}b_m
 - \lf(\fr1{c^2}(bn^k-nb^k)+gq n^k\rg)nn_m
  - g\lf(\fr1{c^2}\Upsilon^k+gq n^k\rg)qb_m.
  $$
Canceling similar terms leaves us with  merely
\be
p^k_m=-\fr1{c^2}B\Upsilon_m^k
-n^2g\fr 1q b^k \fr1{c^2}  \Upsilon_m
+(1-c^2)g\fr bq nn^k \fr1{c^2} \Upsilon_m
-gq nn^k      n_m
 - gq n^k\fr1{c^2}bb_m
  - g^2q^2n^kb_m.
\ee

{%\pgbrk}

\ses

With this formula we find that
\ses

$$
\D{K  m^k}{y^m}-l_mm^k+l^km_m=
 Km^k\fr{q}{2\nu}(1-c^2)g    \fr n{q^3}\Upsilon_m
+
g
\fr1B
\sqrt{\fr q{\nu}}\,\fr cq
\lf(
- b^k \fr{n^2}{c^2}
+(1-c^2)b nn^k \fr1{c^2}
\rg)
\Upsilon_m
 $$

 \ses

 $$
+
g
\fr{gc}{Bq}
\sqrt{\fr q{\nu}}
\lf[
 - \fr{q^2} {c^2}bb_m
  - gq^3b_m
+(b^2+gbq)\fr{1-c^2}{c^2}bb_m
+q^2\fr{1-c^2}{c^2}bb_m
+(b^2+gbq)nn_m
  \rg]
n^k
$$

\ses

\ses

            \ses

$$
=
 Km^k\fr{q}{2\nu}(1-c^2)g    \fr n{q^3}\Upsilon_m
+
g
\fr1B  c
\sqrt{\fr q{\nu}}\,\fr 1q
\lf(
-n^2 b^k \fr1{c^2}
+(1-c^2)b nn^k \fr1{c^2}
\rg)
 \Upsilon_m
 $$

 \ses

\ses

 $$
+
g
\fr{gc}{Bq}
\sqrt{\fr q{\nu}}
(b+gq)
\lf[
b\fr{1-c^2}{c^2}bb_m
-q^2b_m
+bnn_m
  \rg]
n^k,
$$
\ses
or
\ses\\
$$
\D{K  m^k}{y^m}-l_mm^k+l^km_m=
 Km^k\fr{q}{2\nu}(1-c^2)g    \fr n{q^3}\Upsilon_m
$$

\ses

\ses

$$
+
g
\fr1B  c
\sqrt{\fr q{\nu}}\,\fr 1q
n
\lf(
-n b^k \fr1{c^2}
+(1-c^2)b n^k \fr1{c^2}
\rg)
\Upsilon_m
+
g
\fr{gc}{Bq}
\sqrt{\fr q{\nu}}
(b+gq)
n
\Upsilon_m
n^k,
$$
\ses
so that
$$
\D{K  m^k}{y^m}-l_mm^k+l^km_m=
 Km^k\fr{q}{2\nu}(1-c^2)g    \fr n{q^3}(bn_m-nb_m)
+
g
\fr1B
\fr 1q
n
Km^k
 (bn_m-nb_m).
 $$
 \ses
Using here the equality
$$
bn_m-nb_m=cB\fr1K
\sqrt{\fr{q}{\nu}}
m_m
$$
leads to
$$
\D{K  m^k}{y^m}-l_mm^k+l^km_m=
g m^km_m
 \lf[
 \fr{q}{2\nu}(1-c^2)   \fr n{q^3}
B
+
\fr 1q
n
\rg]
c
\sqrt{\fr{q}{\nu}}.
 $$
 \ses
The eventual result reads
\be
\D{K  m^k}{y^m}= l_mm^k -l^km_m+
gc
m^km_m\fr1{2X}
\fr{n}{\sqrt{q\nu}}.
\ee

\ses

{%\pgbrk}

Next, with (B.38)
we find that
$$
\D{\fr1K T}{y^i}=
l_i
\fr1{c}
\sqrt{\fr{\nu} q}\fr{1}{B}
-
\fr1{2c}
\sqrt{\fr q{\nu} }
 (1-c^2)g\fr n{q^3}\Upsilon_i\fr{K}{B}
$$

\ses

\ses

$$
-\fr1{c}
\sqrt{\fr{\nu} q}\fr{K}{B^2}
\lf[
2bb_i + 2
\lf(\fr{1-c^2}{c^2}bb_i+nn_i\rg)
+gb_iq
+\fr1qg
\lf(\fr{1-c^2}{c^2}b^2b_i+bnn_i\rg)\rg],
$$
\ses
or
\be
\D{\fr1K T}{y^i}=
-l_i
\fr1{c}
\sqrt{\fr{\nu} q}\fr{1}{B}
-
\fr1{2c}
\sqrt{\fr q{\nu} }
 (1-c^2)g\fr n{q^3}\Upsilon_i\fr{K}{B}
-\fr1{c} q\nu
\fr 1{\sqrt{\nu q}}
\fr{K}{B^2}
\fr1{q^2}ng
\Upsilon_i.
\ee
Therefore, from  (B.39) it follows that
$$
\D{  m_m}{y^i}=
\Biggl[
-l_i
\fr1{c}
\sqrt{\fr{\nu} q}\fr{1}{B}
-\fr1{c}
\lf[ q\nu+\fr12(1-c^2)B\rg]
\fr 1{\sqrt{\nu q}}
\fr{K}{B^2}
\fr1{q^2}ng
\Upsilon_i
\Biggr]
\Upsilon_m
+
 \fr1{c}
\sqrt{\fr{\nu} q}\fr{K}{B} (b_in_m-b_mn_i)
 $$

 \ses

 \ses

$$
=
-\fr1{c}
\lf[ q\nu+\fr12(1-c^2)B\rg]
\fr 1{\sqrt{\nu q}}
\fr{K}{B^2}
\fr1{q^2}ng
(bn_i-nb_i)
 (bn_m-nb_m)
+
 \fr1{c}
\sqrt{\fr{\nu} q}\fr{K}{B^2}
Y_{im},
$$
where
$$
Y_{im}=B (b_in_m-n_ib_m)-\lf(\fr1{c^2} bb_i+nn_i+gqb_i\rg)(bn_m-nb_m)
=
-\fr BKl_m(bn_i-nb_i).
$$
\ses\\
Using also
$
bn_m-nb_m=c(B/K) \sqrt{{q}/{\nu}} \,m_m,
$
we arrive at the result
\ses\\
\be
K\D{  m_m}{y^i}= -gcm_mm_i\fr1{2X}\fr n{\sqrt{q\nu }}
-
l_mm_i.
\ee

\ses

Differentiating the tensor
$h_{ij}=g_{ij}-l_il_j$ leads to
$$
K\D{h_{ij}}{y^k}=2A_{ijk}-l_ih_{jk}-l_jh_{ik},
$$
where
$A_{ijk}=KC_{ijk}$.
Applying  here (B.57)
yields the representation
$$
A_{ijk}=
-gc\fr1{2X}\fr n{\sqrt{q\nu }}
m_im_jm_k.
$$
We may write
$$
A_{ijk}=
\fr1{A_hA^h}A_iA_jA_k
\quad \text{whenever} ~ g\ne0 ~~ \text{and} ~ n\ne0
.
$$

{%\pgbrk}

We shall use the derivative
\be
\partial^*_n=g_n\D{}g ~ ~ ~ \text{with} ~  ~ g_n=\D g{x^n}.
\ee

When the function $K$ given by (B.18)  is differentiated
with respect to $x^n$, we obtain
\be
\D K{x^n}=\partial^*_n K
+
\fr KBgqy^j\nabla_n b_j
+
a^k{}_{nj}y^jl_k.
\ee
\ses
Let us verify this formula.
With
(B.18)-(B.21)
we find
$$
\D {K^2}{x^n}=\partial^*_n K^2
+
\lf[\partial_n S^2+g(b\partial_n q+q\partial_n b)
-g
(b\partial_n q-q\partial_n b)
\rg]
J^2
=\partial^*_n K^2
+
(\partial_n S^2+2gq\partial_n b)J^2
$$
($\partial_n$ means $\partial/\partial x^n$),
so
that
\ses\\
$$
\D K{x^n}=\partial^*_n K
+
\fr12\fr KBy^iy^j\partial_n a_{ij}
+
\fr KBgq\partial_n b.
$$
We obtain directly that
$$
\fr12\fr KBy^iy^j\partial_n a_{ij}-a^k{}_{nj}y^jl_k
=
\fr12\fr KBy^j\lf[y^i\partial_n a_{ij}
-(y^k+gqb^k)
(\partial_n a_{kj}+\partial_j a_{kn}-\partial_k a_{nj})
\rg],
$$
\ses
or
$$
\fr12\fr KBy^iy^j\partial_n a_{ij}-a^k{}_{nj}y^jl_k
=
-\fr KBy^j
gqb_k
a^k{}_{nj}.
$$
\ses
The formula (B.59) is valid.

\ses

The equality
$
\partial {K^2} /\partial g={\bar M}   K^2
$
holds with
\be
{\bar M}=
 - \fr1{h^3}f+
\fr12\fr{G}{hB}   q^2+  \frac1{h^2B} b q.
\ee
In obtaining this formula we have used the derivatives
$$
  \D hg= -\fr14 G, \quad
    \D Gg= \fr1{h^3}, \quad \D{\lf(\fr Gh\rg)} g =\fr1{h^4} \lf(1+\fr{g^2}4 \rg),
\quad
\D fg= -\fr1{2h} +  \fr{b}B\Bigl(\fr14G  q +   \fr1{2h}b\Bigr).
$$
\ses
Therefore,
\be
\partial^*_nK=
{\bar M}   K^2\D g{x^n}.
\ee
\ses
It follows that
$$
\D {\bar M}{y^h}=\fr{2b^4}{B^2}\D{\fr bq}{y^m}=\fr{4  q^2X}{gBK} A_h
$$
and
$$
\partial^*_nl_m= \D{ \lf(\partial^*_nK\rg)}{y^m}.
$$

All the right-hand parts in the equalities
(B.59)-(B.61)
are regular
of class  $C^{\infty}$ globally regarding the $y$-dependence.

{%\pgbrk}

\ses

To elucidate the global properties of the angle function
$\theta=\theta(x,y)$ in tangent spaces,
we note that  at any point $x\in M$ the space $T_xM\backslash 0$
can conveniently be covered by the atlas
\be
{\cal A} ={\cal C}_1 \cup {\cal C}_2 \cup {\cal C}_3 \cup {\cal C}_4
\ee
 of  the four charts
 \be
 {\cal C}_1\!=\{y\in {\cal C}_1: b>0\}, ~
 {\cal C}_2\!=\{y\in {\cal C}_2: n>0\}, ~
 {\cal C}_3\!=\{y\in {\cal C}_3: n<0\}, ~
 {\cal C}_4\!=\{y\in {\cal C}_4: b<0\}.
  \ee
The north pole and the south pole of the indicatrix (B.25) belong to the charts
${\cal C}_1$ and ${\cal C}_4$, respectively.
In the region
$
{\cal C}_1 \cup {\cal C}_4
$
\ses
the variable
\be
\wt w=\fr n{ b}
\ee
can naturally be used
to introduce the function $\wt\theta(x,\wt w)$ defined by
$\theta=\wt\theta$.
We have
\be
\D{ \wt w}{y^n}=\fr{ b n_n-n b_n}{ b^2}.
\ee
\ses
Comparing this with (B.39) and (B.43),
we conclude
that if
$
 y\in{ \cal C}_1 \cup {\cal C}_4,
 $
 then
\be
\D{\wt \theta}{\wt w}
=
\fr1c  \fr{b^2}B
\sqrt{  \fr{\nu}q   }.
\ee
In alternative regions  the variable
\be
t =-\fr { b}q, ~ ~ \text{if} ~ ~ y\in{ \cal C}_2;  \quad
t =\fr { b}q, ~ ~ \text{if} ~ ~ y\in{ \cal C}_3,
\ee
is well adaptable.
Since
\be
\D{t}{y^m}=    \fr{| n|}{q^3}(bn_m-nb_m)
\ee
(see (B.11)),
we obtain
that if
$
 y\in{ \cal C}_2 \cup {\cal C}_3,
 $
 then
 \be
\D{\wh \theta}{t }
=
\fr1c  \fr{q^3}{|n|B}
\sqrt{  \fr{\nu}q   } ,
\ee
when
using the function
 $\wh\theta(x,t)$ defined by
$\theta=\wh\theta$.
\ses
The positivity
\be
\D{\wt \theta}{\wt w}  >0,
\quad
\D{\wh \theta}{t }>0
\ee
holds.

{%\pgbrk}

It can readily be seen that in intersections of the charts introduced
the variables $\wt w$ and $t$ are expressible one through another
in the $C^{\infty}$-smooth way,
so that  the atlas $\cA$ introduced in (B.62) is smooth of the class
$C^{\infty}$.
The right-hand parts of the angle derivatives (B.66) and (B.69)
are also of this class with respect to $y$.
Noting also the positivity (B.70), we are entitled to conclude that
 when point moves along the indicatrix from the north pole  in the right
direction
the value $\theta$ of the point
increases monotonically
\be
0\le\theta< \theta^{\text{max}},
\ee
where
\be
\theta^{\text{max}}=
2(\theta^{I}+\theta^{II})
\ee
with
\be
\theta^{I}=
\int_0^{\infty}
\fr1c  \fr{b^2}B
\sqrt{  \fr{\nu}q   }d\wt w,
  \qquad
\theta^{II}=
\int_{-\infty}^0
\fr1c  \fr{b^2}B
\sqrt{  \fr{\nu}q   }d\wt w.
\ee
Reminding the expressions  of the functions $B$ and $\nu$
(see (B.14) and (B.26)),
we can write these integrals explicitly as follows:
\be
\theta^{I}=
\fr1c
\int_0^{\infty}
\fr{\sqrt{1+(1-c^2)g\fr1w}}{     1+gw+w^2        }
d\wt w,
 \qquad
\theta^{II}=
\fr1c
\int_0^{\infty}
\fr{\sqrt{1-(1-c^2)g\fr1w}}{     1-gw+w^2        }
d\wt w,
 \ee
where
\be
w=\sqrt{\wt w^2+\fr{1-c^2}{c^2}}\equiv \fr q{|b|}.
\ee
The value $\theta=0$ corresponds to the north pole of indicatrix, and
the
value
$\theta=\theta^{I}+\theta^{II}$ to the south pole.
The indicatrix intersects the direction with $b=0, n>0$
at $\theta=\theta^{I}$,
so that the angle measure of the upper chart
${ \cal C}_1 $
 equals $2\theta^{I}$, respectively
 $2\theta^{II}$ for the chart
${ \cal C}_4 $.
The value $\theta^{\text{max}} =2(\theta^{I}+\theta^{II})$
is the total length of the indicatrix.
The vertical straight angle going from the north to the south
costs
$
\theta^{I}+\theta^{II}.
$

In general,
$\theta^{\text{max}}=\theta^{\text{max}}(x).
$

In the Riemannian limit we obtain
$
\theta^{I}=\theta^{II}=\pi/2
$
(put $g=0$ and $c=1$ in (B.74)).

{%\pgbrk}

In the region
${ \cal C}_1$
we can interpret
the $\wt\theta$ to be of the functional dependence
\be
\wt\theta=\wt\Theta\bigl(g(x), c(x),\wt w\bigr),
\ee
obtaining from  (B.66)
\be
\wt\Theta
=
\fr1c
\int
\fr{\sqrt{1+(1-c^2)g\fr1w}}
{     1+gw+w^2        }
d\wt w.
\ee
The integration constant should be specified by means of the condition
\be
\wt\theta_{\bigl|\wt w=0\bigr.}=0
\ee
to be in agreement with $\theta(x,b(x))=0$.
Differentiating the equality
 $\theta(x,y)=\wt\Theta\bigl(g(x), c(x),\wt w\bigr)$
with respect to $x^n$ yields
\be
\D{\theta}{x^n}
=
\D{\wt\Theta}g
\D g{x^n}
+\D{\wt\Theta}c
\D c{x^n}
+\D{\wt\Theta}{\wt w}
\D {\wt w}{x^n}.
\ee
\ses
From (B.77) it follows that
\be
\D{\wt \Theta}g
=
\fr1c
\int
\lf(
-  \fr{w\,\sqrt{1+(1-c^2)g\fr1w}}
{  (1+gw+w^2  )^2}
+
 \fr{(1-c^2)}{2w(1+gw+w^2 )}
\fr1{\sqrt{1+(1-c^2)g\fr1w}}
\rg)
d \wt w
\ee
and
\be
\D{\wt \Theta}c
=
-\fr1{c}\wt\Theta
-
g\int \fr1w
\fr1{(     1+gw+w^2 )  \sqrt{1+(1-c^2)g\fr1w}
}
d\wt w.
\ee
Using (B.76), (B.66), and (B.39), we can write
\be
\D{\theta}{x^n}
=
\D{\wt\Theta}g
\D g{x^n}
+\D{\wt\Theta}c
\D c{x^n}
+
\fr1{K^2}Tb^2
          \lf(
 \D{\wt w}{x^n}
- a^k{}_{nj}y^j
(bn_k-nb_k)
\fr1{b^2}
\rg)
+
  a^k{}_{nj}y^j\D{\theta}{y^k}.
\ee
\ses
Taking into account equality
$$
 \D{\wt w}{x^n}=
 y^k\fr1{b^2}
 (b\nabla_nn_k-n\nabla_nb_k)
+
\fr1{b^2}
 (bn_k-nb_k)
a^k{}_{nj}y^j,
 $$
we obtain
$$
 \D{\wt w}{x^n}
- a^k{}_{nj}y^j
(bn_k-nb_k)
\fr1{b^2}
=
 y^k\fr1{b^2}
 (b\nabla_nn_k-n\nabla_nb_k),
 $$
which yields
$$
b^2
\lf( \D{\wt w}{x^n} - a^k{}_{nj}y^j (bn_k-nb_k) \fr1{b^2}\rg)
=
-S^2cn^h\nabla_n  \wt b_h-nb\fr1c\D c{x^n}.
 $$
\ses
Inserting this result in (B.82),  we are coming to
\be
\D{\theta}{x^n}
=
\D{\wt\Theta}g
\D g{x^n}
+\D{\wt\Theta}c
\D c{x^n}
-
\fr1{K^2} T\lf(S^2cn^h\nabla_n  \wt b_h+bn\fr1c\D c{x^n}\rg)
+
  a^k{}_{nj}y^j\D{\theta}{y^k}.
\ee
\ses
Here, all the terms are smooth
  of class $C^{\infty}$
regarding the $y$-dependence.
Similar representations for the derivative $\partial{\theta}/\partial{x^n}$
can be obtained in the regions
${ \cal C}_2,{ \cal C}_3,{ \cal C}_4$.

{%\pgbrk}

\ses

Henceforth, we assume
\be
c=const.
\ee
Let us propose the expansion
\be
N^k_n=
\Biggl[
\lf(
\fr1{c^2}b
-
Uc^2B_1
\rg)
n^k
+
n
\lf(U-\fr1{c^2}\rg)
b^k
\Biggr]
p_n
-l^k\partial^*_nK-Km^k\partial^*_n\theta
-
  a^k{}_{nj}y^j,
\ee
where
$
B_1=(1/c^2)b+gq
$
(in agreement with  (B.14))
and
\be
p_n=cn^h\nabla_n  \wt b_h
\equiv
cn^h
\lf(\D{\wt b_h}{x^n}-a^k{}_{nh}\wt b_k\rg).
\ee

\ses

By contracting the coefficients  $N^k_n$ written in (B.85) we obtain
\be
N^k_ny_k=-K\partial^*_nK-gnqp_n
\fr{K^2}B
  -  a^k{}_{nj}y^jl_k.
  \ee
\ses
On the other hand,
\be
\D K{x^n}=\partial^*_n K
+
\fr KBgnqp_n
+
a^k{}_{nj}y^jl_k
\ee
(see (B.59)).
Comparing this equality with (B.87) just yields
\ses
\be
d_nK=0.
\ee

{%\pgbrk}

If we take in (B.85) the function $U$ to be
\be
U
=C_1\fr1{c}
\sqrt{\fr q{\nu}}, \qquad   C_1=-\fr1c,
\ee
we obtain
\be
\D{\theta}{x^n}+N^k_n\D{\theta}{y^n}
=-C_1p_n.
\ee
Indeed,
$$
N^k_n \D{\theta}{y^k}=
\fr1K N^k_n m_k
=
-\partial^*_n\theta
+
\fr1{K^2} T N^k_n (bn_k-nb_k)
-
  a^k{}_{nj}y^j\D{\theta}{y^k},
$$
\ses
or
$$
N^k_n \D{\theta}{y^k}=
-\partial^*_n\theta
+
\fr1{K^2} T(S^2-Uc^2B)p_n
 -
  a^k{}_{nj}y^j\D{\theta}{y^k},
$$
which can also be written as follows:
\ses
\be
N^k_n \D{\theta}{y^k}=
-\partial^*_n\theta
+
\fr1{K^2} TS^2p_n
-C_1p_n
 -
  a^k{}_{nj}y^j\D{\theta}{y^k}.
  \ee
\ses
On the other hand,
\be
\D{\theta}{x^n}
=
\partial^*_n\theta
-
\fr1{K^2} TS^2p_n
+
  a^k{}_{nj}y^j\D{\theta}{y^k},
\ee
where  (B.83) has been used.
By adding this equality to (B.92) we just conclude that
the choice (B.85) made for the coefficients $N^k_n$ entails
the angle-preserving
property $d_n\theta=k_n$ (indicated in (2.8))
with the choice
 $k_n=-C_1p_n$.

\ses

\ses
{%\pgbrk}

Let us insert  (B.88) and (B.93) in (B.85), which yields
\ses\\
$$
N^k_n=
\Biggl[
\lf(
\fr1{c^2}b
-
Uc^2B_1
\rg)
n^k
+
n
\lf(U-\fr1{c^2}\rg)
b^k
\Biggr]
p_n
+
l^k
\lf[
 gnqp_n   \fr{K}B
+
  a^h{}_{nj}y^jl_h
  \rg]
$$

\ses

$$
-
Km^k
\lf[
\fr1{K^2} TS^2p_n
-
  a^h{}_{nj}l^jm_h
\rg]
-l^k\partial_nK-Km^k\partial_n\theta
-
  a^k{}_{nj}y^j.
  $$
\ses
Noting  the identity
$l^kl_h+m^km_h=\de^k_h$
and
the representation (B.90) of the function $U$,
\ses
we get
$$
N^k_n=
\fr1{c^2}(bn^k-nb^k)
p_n
+
l^k
 gnqp_n   \fr{K}B
-
Km^k
\fr1{K^2} TS^2p_n
-l^k\partial_nK-Km^k(\partial_n\theta+C_1p_n).
  $$
\ses
Using
here the formula (B.40) which describes the structure of the vector $m^k$,
we obtain
$$
N^k_n=
\lf[
B\fr1{c^2}(bn^k-nb^k)
+
y^k gnq
\rg]
\fr1B
p_n
-
\lf(
\fr1{c^2}
(bn^k-nb^k)
+
gq
 n^k
\rg)
\fr{1}{B}
S^2p_n
$$

\ses

$$
-l^k\partial_nK-Km^k(\partial_n\theta+C_1p_n)
$$

\ses

\ses

\ses

$$
=
\lf[
gbq\fr1{c^2}(bn^k-nb^k)
+
y^k gnq
\rg]
\fr1B
p_n
-
gq
 n^k
\fr{1}{B}
S^2p_n
-l^k\partial_nK-Km^k(\partial_n\theta+C_1p_n)
$$

\ses

\ses

\ses

$$
=
\lf[
gbq\fr1{c^2}bn^k
+
nn^k gnq
\rg]
\fr1B
p_n
-
gq
 n^k
\fr{1}{B}
S^2p_n
-l^k\partial_nK-Km^k(\partial_n\theta+C_1p_n),
$$
\ses
so that
the final representation is merely
\be
N^k_n=
-l^k\D K{x^n}-Km^kP_n
\ee
with
\be
P_n=\D{\theta}{x^n}+C_1p_n.
\ee
We have arrived at the coefficients $N^k_n$ which are tantamount to the coefficients
announced in (2.6).

{%\pgbrk}

With the coefficients $N^k_n$ indicated in
(B.85), we obtain the values
\ses\\
\be
d_ib=nUc^2p_i-\fr1K b \partial^*_iK
+
c\sqrt{\fr q{\nu}}
n
  \partial^*_i\theta
\ee
\ses
 and
\be
d_i n=-   Uc^2B_1p_i
-\fr1K n \partial^*_iK
-
c\sqrt{\fr q{\nu}}
B_1
  \partial^*_i\theta.
\ee

Also, using the variable
$$
q=\sqrt{\fr{1-c^2}{c^2}b^2+n^2},
$$
we get
$$
qd_i q =  \fr{1-c^2}{c^2}b \lf[
nUc^2p_i-\fr1K b \partial^*_iK
+
c\sqrt{\fr q{\nu}}
n
  \partial^*_i\theta
    \rg]
+n
\lf[
-   Uc^2B_1p_i
-\fr1K n \partial^*_iK
-
c\sqrt{\fr q{\nu}}
B_1
  \partial^*_i\theta  \rg].
$$

For the square
$$
S^2=b^2+q^2
$$
\ses
the equality
$$
\fr12d_i S^2 =  \fr{1}{c^2}b \lf[
nUc^2p_i-\fr1K b \partial^*_iK
+
c\sqrt{\fr q{\nu}}
n
  \partial^*_i\theta
    \rg]
+n
\lf[
-   Uc^2B_1p_i
-\fr1K n \partial^*_iK
-
c\sqrt{\fr q{\nu}}
B_1
  \partial^*_i\theta  \rg]
$$
\ses
is obtained,
so that
\be
\fr12d_i S^2 =
-   Uc^2 gnq p_i
-S^2\fr1K  \partial^*_iK
-
c\sqrt{\fr q{\nu}}
gnq
  \partial^*_i\theta
=
-S^2\fr1K  \partial^*_iK
-
c\sqrt{\fr q{\nu}}
gnq
P_i
\ee
\ses
and
\be
qd_i q =
-q^2\fr1K  \partial^*_iK
-
c\sqrt{\fr q{\nu}}
(b+gq)
nP_i.
\ee

It can readily be verified that
$$
d_i\fr nq= -\fr{1-c^2}{c^2q^3}bB
c\sqrt{\fr q{\nu}} P_i.
$$

{%\pgbrk}

We directly obtain the equalities
$$
d_i(bq)=
\fr bq
 \fr{1}{c^2}b \lf[
nUc^2p_i-\fr1K b \partial^*_iK
+
c\sqrt{\fr q{\nu}}
n
  \partial^*_i\theta
    \rg]
-
\fr bq
b \lf[
nUc^2p_i-\fr1K b \partial^*_iK
+
c\sqrt{\fr q{\nu}}
n
  \partial^*_i\theta
    \rg]
$$

\ses

$$
+
\fr bq
n
\lf[
-   Uc^2B_1p_i
-\fr1K n \partial^*_iK
-
c\sqrt{\fr q{\nu}}
B_1  \partial^*_i\theta  \rg]
+qnUc^2p_i-
q\fr1K b \partial^*_iK
+
qc\sqrt{\fr q{\nu}}
n
  \partial^*_i\theta,
  $$
\ses
or
$$
d_i(bq)=
-
\fr bq
b \lf[
nUc^2p_i
+
c\sqrt{\fr q{\nu}}
n
  \partial^*_i\theta
    \rg]
+
\fr bq
n
\lf[
-   Uc^2gqp_i
-
c\sqrt{\fr q{\nu}}
gq   \partial^*_i\theta  \rg]
 $$

\ses

\ses

$$
+qnUc^2p_i-
2bq\fr1K  \partial^*_iK
+
qc\sqrt{\fr q{\nu}}
n
  \partial^*_i\theta,
  $$
\ses
and
\be
d_i B=
-
g\fr 1q
nB
Uc^2p_i
-
g\fr 1q
B
c\sqrt{\fr q{\nu}}
n
  \partial^*_i\theta
-
2B
\fr1K  \partial^*_iK
+bqg_i,
\ee

{%\pgbrk}

\ses

\nin
together with
$$
d_i \fr bq=
-\fr b{q^3}
 \fr{1}{c^2}b \lf[
nUc^2p_i-\fr1K b \partial^*_iK
+
c\sqrt{\fr q{\nu}}
n
  \partial^*_i\theta
    \rg]
+
\fr b{q^3}
b \lf[
nUc^2p_i-\fr1K b \partial^*_iK
+
c\sqrt{\fr q{\nu}}
n
  \partial^*_i\theta
    \rg]
$$

\ses

$$
-
\fr b{q^3}
n
\lf[
-   Uc^2 B_1 p_i
-\fr1K n \partial^*_iK
-
c\sqrt{\fr q{\nu}}
B_1  \partial^*_i\theta  \rg]
+\fr1qnUc^2p_i-
\fr1q\fr1K b \partial^*_iK
+
\fr1qc\sqrt{\fr q{\nu}}
n
  \partial^*_i\theta,
  $$
\ses
which is
$$
d_i \fr bq=
\fr b{q^3}
 \fr{1}{c^2}b
\fr1K b \partial^*_iK
+
\fr b{q^3}
b \lf[
-\fr1K b \partial^*_iK
+
c\sqrt{\fr q{\nu}}
n
  \partial^*_i\theta
    \rg]
$$

\ses

$$
+
\fr b{q^3}
n
\lf[
\fr1K n \partial^*_iK
+
c\sqrt{\fr q{\nu}}
gq
  \partial^*_i\theta  \rg]
+\fr1{q^3}nB Uc^2p_i-
\fr1q\fr1K b \partial^*_iK
+
\fr1qc\sqrt{\fr q{\nu}}
n \partial^*_i\theta,
  $$
\ses
coming to
\be
d_i \fr bq=
\fr n{q^3}B
\lf(Uc^2p_i
+
c\sqrt{\fr q{\nu}}
  \partial^*_i\theta
  \rg)
=  \fr n{q^3}B
c\sqrt{\fr q{\nu}}\lf(C_1p_i
+
  \partial^*_i\theta
  \rg).
 \ee

\ses

\ses

Applying the operator $d_i$ to the scalar $T$ introduced by (B.38)
yields
$$
d_iT=T\lf(\fr12\fr q{\nu}d_i\fr{\nu}q-\fr1Bd_iB\rg)
=T\fr12\fr q{\nu}(1-c^2)\lf(g_i\fr bq +
g
\fr n{q^3}B
c\sqrt{\fr q{\nu}}\lf(C_1p_i
+
  \partial^*_i\theta
  \rg)
  \rg)
$$

\ses

$$
-
T\fr1B\lf(
-
g\fr 1q
nB
Uc^2p_i
-
g\fr 1q
B
c\sqrt{\fr q{\nu}}
n
  \partial^*_i\theta
-
2B
\fr1K  \partial^*_iK
+bqg_i
 \rg),
$$
\ses
or
\ses\\
\be
\fr1Td_iT=\fr {bq}{2B}\lf(\fr1X-4\rg)g_i
+
2
\fr1K  \partial^*_iK
+
g
\fr n{2q}\fr1X
\lf(
Uc^2
p_i
+
c\sqrt{\fr q{\nu}}
  \partial^*_i\theta
\rg).
\ee

\ses

We also obtain

$$
d_h
\Biggl[
\lf(\fr nq\rg)^2 \fr q{\nu}
\Biggr]
=
-2
\fr nq \fr q{\nu}\fr{1-c^2}{c^2q^3}bB
c\sqrt{\fr q{\nu}}
P_h
-
\lf(\fr nq\rg)^2 \lf( \fr q {\nu}\rg)^2
d_h\fr{\nu}q,
$$
\ses
or
$$
d_h
\Biggl[
\lf(\fr nq\rg)^2 \fr q{\nu}
\Biggr]
=
-2
\fr n{\nu}\fr{1-c^2}{c^2q^3}bB
c\sqrt{\fr q{\nu}}
P_h
-
 \lf( \fr n {\nu}\rg)^2
(1-c^2)
\lf(g_h\fr bq +
g
\fr n{q^3}B
c\sqrt{\fr q{\nu}}
P_h
  \rg),
  $$
\ses
and
$$
d_h
\fr1X
=
2
\fr n{\nu}\fr{1-c^2}{q^3}bB
c\sqrt{\fr q{\nu}}
P_h
+
 \lf( \fr n {\nu}\rg)^2
c^2(1-c^2)
\lf(g_h\fr bq +
g
\fr n{q^3}B
c\sqrt{\fr q{\nu}}
P_h
  \rg).
  $$

  {%\pgbrk}

\ses

Next,
starting from
 (B.40),
it follows that
$$
d_i(Km^n)=
-
Km^n \fr12\fr q{\nu}d_i\fr{\nu}q
+c
\sqrt{\fr q{\nu}}
d_i
\lf(
\fr1{c^2}
(bn^n-nb^n)
+
gq
 n^n
\rg),
$$
\ses
or
$$
d_i(Km^n)=
-
Km^n
\fr12\fr q{\nu}(1-c^2)\lf(g_i\fr bq +
g
\fr n{q^3}B
c\sqrt{\fr q{\nu}}\lf(C_1p_i
+
  \partial^*_i\theta
  \rg)
  \rg)
  $$

\ses

$$
+
c
\sqrt{\fr q{\nu}}
\fr1{c^2} n^n
\lf[
nUc^2p_i-\fr1K b \partial^*_iK
+
c\sqrt{\fr q{\nu}}
n
  \partial^*_i\theta
  \rg]
  $$

\ses

\ses

$$
-
c
\sqrt{\fr q{\nu}}
\fr1{c^2} b^n
\lf[
-   Uc^2\lf( \fr1{c^2}b+gq\rg)p_i
-\fr1K n \partial^*_iK
-
c\sqrt{\fr q{\nu}}
\lf( \fr1{c^2}b+gq\rg)
  \partial^*_i\theta
  \rg]
$$

\ses

$$
+
c
\sqrt{\fr q{\nu}}
\fr1{c^2} (b\nabla_in^n-n\nabla_ib^n)
+
c
\sqrt{\fr q{\nu}}
\lf(
gq \nabla_in^n
+g_iq n^n
\rg)
$$

\ses

\ses

$$
+
c
\sqrt{\fr q{\nu}}
g n^n
\fr1q
  \fr{1-c^2}{c^2}b \lf[
nUc^2p_i-\fr1K b \partial^*_iK
+
c\sqrt{\fr q{\nu}}
n
  \partial^*_i\theta
    \rg]
$$

\ses

$$
+
c
\sqrt{\fr q{\nu}}
g n^n
\fr1q
n
\lf[
-   Uc^2\lf( \fr1{c^2}b+gq\rg)p_i
-\fr1K n \partial^*_iK
-
c\sqrt{\fr q{\nu}}
\lf( \fr1{c^2}b+gq\rg)
  \partial^*_i\theta  \rg],
$$

\ses

\ses

\nin
which is transformed to
\ses\\
$$
d_i(Km^n)=
-
Km^n
\fr12\fr q{\nu}(1-c^2)\lf(g_i\fr bq +
g
\fr n{q^3}B
c\sqrt{\fr q{\nu}}\lf(C_1p_i
+
  \partial^*_i\theta
  \rg)
  \rg)
  $$

\ses

$$
+
c
\sqrt{\fr q{\nu}}
\fr1{c^2} n^n
\lf[
-\fr1K b \partial^*_iK
+
c\sqrt{\fr q{\nu}}
n
  \partial^*_i\theta
  \rg]
  $$

\ses

\ses

$$
-
c
\sqrt{\fr q{\nu}}
\fr1{c^2} b^n
\lf[
-\fr1K n \partial^*_iK
-
c\sqrt{\fr q{\nu}}
\lf( \fr1{c^2}b+gq\rg)
  \partial^*_i\theta
  \rg]
$$

\ses

$$
+
c
\sqrt{\fr q{\nu}}
\fr1{c^2} (b\nabla_in^n-n\nabla_ib^n)
+
c
\sqrt{\fr q{\nu}}
\lf(
gq \nabla_in^n
+g_iq n^n
\rg)
+
c
\fr q{\nu}
g n^n
  \fr{1-c^2}{qc^2}b
c
n
  \partial^*_i\theta
$$

\ses

$$
+
c
\sqrt{\fr q{\nu}}
g n^n
\fr1q
\lf[
-\fr1K q^2 \partial^*_iK
-
nc\sqrt{\fr q{\nu}}
\lf( \fr1{c^2}b+gq\rg)
  \partial^*_i\theta  \rg]
$$

\ses

\ses

\ses

$$
+
Uc^2
c \sqrt{\fr q{\nu}}
\fr1{c^2}
\lf[
\fr nq (q-gc^2(b+gq))    n^n
+b^n
\lf( \fr1{c^2}b+gq\rg)
 \rg]
p_i.
$$
\ses
Notice that
$$
b\nabla_in^n-n\nabla_ib^n
=
-y^np_i,  \qquad
\nabla_in^n=-\fr1{c^2}b^np_i.
$$

{%\pgbrk}

\nin
In this way  we straightforwardly obtain the representation
\ses
$$
d_i(Km^n)=
-
Km^n
\fr12\fr q{\nu}(1-c^2)\lf(g_i\fr bq +
g
\fr n{q^3}B
c\sqrt{\fr q{\nu}}\,P_i
  \rg)
-
\fr1c\sqrt{\fr q{\nu}}
\,
y^np_i
+
c
\sqrt{\fr q{\nu}}
g_iq n^n
  $$

\ses

\ses

\be
-
c
\sqrt{\fr q{\nu}}\,
\fr1{c^2}gqb^np_i
-
m^n  \partial^*_iK
+\fr q{\nu}
\lf[y^n
-\fr nq gc^2(b+gq)    n^n
+
gqb^n
 \rg]
P_i.
\ee

{%\pgbrk}

\ses

We can write the coefficients  (B.85) in terms of the vector $m^k$
\be
N^k_n=-Km^k(C_1p_n+\partial^*_n\theta)
+
\fr1{c^2}
(bn^k-nb^k)
p_n
-l^k\partial^*_nK
-
  a^k{}_{nj}y^j
\ee
and after that represent them
in the form
\be
N^k_n=
-
\Bigl(N_{\{l\}}l^k+N_{\{m\}}m^k\Bigr)
Kp_n
-
l^k\partial^*_nK
-Km^k\partial^*_n\theta
-
  a^k{}_{nj}y^j
\ee
\ses
with
\be
N_{\{l\}}
=
g\fr{ nq}{B}, \qquad
N_{\{m\}}=
C_1-      \fr1c\sqrt{\fr {\nu}q}         \fr{S^2}{B}.
\ee

{%\pgbrk}

\ses

Differentiating the coefficients $N^k_n$ given by (B.85)
leads to the representation
\ses\\
$$
N^k_{nm}=
\Biggl[
b_m
\lf(
\fr1{c^2}-
U
\fr{\nu}q
\rg)
n^k
+
n_m
\lf(U-\fr1{c^2}\rg)
b^k
\Biggr]
p_n
-
gUc^2
\fr nq
n_m
n^k
p_n
$$

\ses
\ses

\ses

$$
+
\fr { n}{2q^2}
\fr {U}{\nu}
g
(1-c^2)
z^k
(bn_m-nb_m)
p_n
-
a^k{}_{nm}
$$

\ses

\ses

\ses

$$
-l^k\partial^*_nl_m
+ \lf(-l_mm^k+l^km_m -
gc
m^km_m\fr1{2X}
\fr{n}{\sqrt{q\nu}}
\rg)
\partial^*_n\theta
-m^k\partial^*_nm_m.
$$
\ses
Here,
$$
z^k=     nb^k-(b+gc^2q)n^k =-\sqrt{\fr{\nu}q} \,cKm^k,  \qquad
 U
=C_1\fr1{c}
\sqrt{\fr q{\nu}},
$$
\ses
so that we can write the above coefficients as follows:
$$
N^k_{nm}=
\Biggl[
b_m
\lf(
\fr1{c^2}-
C_1\fr1{c}
\sqrt{\fr {\nu}q}\,
\rg)
n^k
+
n_m
\lf(C_1\fr1{c}
\sqrt{\fr q{\nu}}-\fr1{c^2}\rg)
b^k
\Biggr]
p_n
-
gC_1\fr1{c}
\sqrt{\fr q{\nu}}c^2
\fr nq
n_m
n^k
p_n
$$

\ses
\ses

\ses

$$
-
\fr { n}{2q^2}
\fr {1}{\nu}
g
(1-c^2)
C_1
Km^k
(bn_m-nb_m)
p_n
-
a^k{}_{nm}
$$

\ses

\ses

\ses

\be
-l^k\partial^*_nl_m
+ \lf(-l_mm^k+l^km_m -
gc
m^km_m\fr1{2X}
\fr{n}{\sqrt{q\nu}}
\rg)
\partial^*_n\theta
-m^k\partial^*_nm_m,
\ee
\ses
or
$$
N^k{}_{nm}=
\lf(-l_mm^k+l^km_m -
gc
m^km_m\fr1{2X}
\fr{n}{\sqrt{q\nu}}
\rg)
P_n
-l^k\partial^*_nl_m
-m^k\partial^*_nm_m
$$

\ses

\be
+
\fr1{c^2}
(b_mn^k-n_mb^k)
p_n-
  a^k{}_{nm}.
  \ee
\ses
This representation can be written in the form similar to (2.11).

{%\pgbrk}

From the previous representation
the  coefficients
$
N^k_{nmi}=\partial{N^k_{nm}}/\partial{y^i}
$
are evaluated to read
\be
N^k_{nmi}=\fr1K Z_nm^km_mm_i
\ee
with
\be
Z_{n}
=
-
\fr 34
g
(1-c^2)^2
B^2
 \fr 1{q^3\nu^3}
(2 b\nu+gc^2n^2)
P_n
+
\fr12
nc\partial^*_n
\lf(g\fr1{X}\fr 1{\sqrt{q\nu }}\rg).
\ee
The last two  formulas  agree completely
 with the general
formula (2.14) indicated in Section 2.

\ses

The following equalities
$$
\D{\fr{\nu}q}{y^i}=(1-c^2)g\fr n{q^3}(nb_i-bn_i), \qquad
\D{\fr n{\nu}}{y^i}=           \fr 1{q\nu^2} \fr1{c^2}(1-c^2)(b+gc^2q)(bn_i-nb_i),
$$

\ses

 $$
 \D{\fr1X}{y^i}=
-
(1-c^2)
(2 b\nu+gc^2n^2)
\fr 1{\nu^2}\fr n{q^3}cB\fr1K \sqrt{\fr{q}{\nu}} m_i,
$$

\ses

$$
\D{\lf(\fr{n^2}{q\nu}\rg)}{y^i}=
\fr 1{c^2}(1-c^2)
\fr{n}{\nu}(bn_i-nb_i)
 \fr 1{q^3\nu}
(2 b\nu+gc^2n^2),
$$
\ses
and
$$
\D{A^hA_h}{y^i}=
\fr{3g^2}{4X}
(1-c^2)^2
(2 b\nu+gc^2n^2)
\fr 1{\nu^3}\fr 1{q^3}B^2\fr1K \fr{cn}{\sqrt{q\nu}} m_i
$$
 are convenient to use in process of derivation of the $Z_n$
indicated above.

{%\pgbrk}

\ses

It is easy again to perform   the evaluation of the curvature tensor
making the choice of the vector $m^k$
in accordance with (B.38)-(B.40).
With arbitrary smooth $g=g(x)$ and $c=c(x)$, the result reads
\be
 M^n{}_{ij}= d_iN^n_j-d_jN^n_i
=Km^nM_{ij}, \qquad
M_{ij}=\D{k_j}{x^i}-\D{k_i}{x^j}
\ee
and
\be
E_k{}^n{}_{ij}=
-
\lf[
 l_km^n -l^nm_k+
gc
m^nm_k\fr1{2X}
\fr{n}{\sqrt{q\nu}}
\rg]
M_{ij}
\equiv
-\D{M^n{}_{ij}}{y^k}
\ee
((B.55) has been applied),
so that the  curvature tensor
\be
\Rho_k{}^n{}_{ij}=E_k{}^n{}_{ij}
-
M^h{}_{ij}C^n{}_{hk}
\ee
is simply
\be
\Rho_k{}^n{}_{ij}=
(l^n m_k-l_k m^n)
M_{ij}.
\ee

Let us find
$$
l_n m_k-l_k m_n=
 T\fr{1}B
 \lf[   ( u_n+gqb_n )  (bn_k-nb_k)- ( u_k+gqb_k )  (bn_n-nb_n)  \rg]
$$

\ses

$$
=
 T\fr{1}B
 \lf[    \wt b^2  (b_nn_k-b_kn_n)  - n^2 (b_kn_n-b_nn_k)
 +gqb(b_n n_k-b_k n_n)
    \rg].
$$
Thus,
\be
l_n m_k-l_k m_n=
 T
(b_n n_k-b_k n_n)
\ee
and
\be
\Rho_{knij}=
 T
(b_n n_k-b_k n_n)
M_{ij},
\ee
so that the factor function $f_1$ indicated in (4.20) is equal to
$cT$.

{%\pgbrk}

\ses

\ses

\setcounter{equation}{0}

\bc{ \bf Appendix C: Randers metric  in two-dimensional case} \ec

\ses

\ses

The Finsler metric function $F$ is now of the form
\be
F=S+b.
\ee
All the formulas (B.1)--(B.11) can be applied.
We obtain
$
\det(g_{ij})=
( F/S)^3
\det(a_{ij}),
$
so that
the
 formulas (B.38)--(B.40) are to be replaced by
\be
T=
\fr1c
\sqrt{\fr{\det(g_{ij})}{\det(a_{mn})}}
=
\fr1c\sqrt{\lf(\fr FS \rg)^3}
\ee
and
\be
  m_i=
   \fr1{F}T(bn_i-nb_i)
   = \fr1{cS}\sqrt{\fr FS} \, (bn_i-nb_i),
\ee
together with
 $m^i=g^{ij}m_j$
 and
 \be
m^i= \fr1{cF}\sqrt{\fr SF}
\lf(   -nb^i  +(b+c^2S)  n^i \rg).
\ee
The partial derivative
$
\partial {F}/\partial{x^n}
$
of the $F$ given by (C.1) is obviously regular of class $C^{\infty}$
with respect to the variable $y$.

Again, we can use the atlas (B.62)--(B.63) and the formulas of the type
(B.64)--(B.70).
In place of (B.77) we obtain the representation
$$
\wt\Theta
=
\int \fr{b^2}{F^2}
Td\wt w
=
\fr1c
\int \fr{1}{\sqrt{1+\sqrt{1+w^2}}\lf(\sqrt{1+w^2}\rg)^3}
d\wt w,
$$
or
\be
\wt\Theta
=
\fr1c
\int \fr{1}{\sqrt{1+\sqrt{\fr1{c^2}+\wt w^2}}\lf(\sqrt{\fr1{c^2}+\wt w^2}\rg)^3}
d\wt w,
\ee
where  the integration constant is assumed to be subjected to the
 condition
$ \wt\Theta_{\bigl|\wt w=0\bigr.}=0
$
to agree with $\theta(x,b(x))=0$.

It follows that
\be
\D{\theta}{x^n}
=
\D{\wt\Theta}c
\D c{x^n}
+\D{\wt\Theta}{\wt w}
\D {\wt w}{x^n}
=
\D{\wt\Theta}c
\D c{x^n}
+
\fr1{F^2}Tb^2
          \lf(
 \D{\wt w}{x^n}
- a^k{}_{nj}y^j
(bn_k-nb_k)
\fr1{b^2}
\rg)
+
  a^k{}_{nj}y^j\D{\theta}{y^k},
\ee
\ses
so that
\be
\D{\theta}{x^n}
=
\D{\wt\Theta}c
\D c{x^n}
-
\fr1{F^2} T\lf(S^2cn^h\nabla_n  \wt b_h+bn\fr1c\D c{x^n}\rg)
+
  a^k{}_{nj}y^j\D{\theta}{y^k}.
\ee
\ses
Here, all the terms are smooth
  of class $C^{\infty}$
regarding the $y$-dependence.
The same conclusion can be arrived at
in the regions
${ \cal C}_2,{ \cal C}_3,{ \cal C}_4$.

The factor-relation (4.19) for the curvature tensor
takes now on
the form
\be
\Rho_{knij}=
\sqrt{\lf(\fr FS \rg)^3}
\bar L_{knij}  \quad \text{with} ~ ~
\bar L_{knij}= a_{nh}\bar L_k{}^h{}_{ij}
\equiv
\bar L_{knij}(x).
\ee

{%\pgbrk}

\vskip 1cm

\def\bibit[#1]#2\par{\rm\noindent\parskip1pt
                     \parbox[t]{.05\textwidth}{\mbox{}\hfill[#1]}\hfill
                     \parbox[t]{.925\textwidth}{\baselineskip11pt#2}\par}

\bc {\bf  References}
\ec

\ses

\bibit[1] H. Rund, \it The Differential Geometry of Finsler
 Spaces, \rm Springer, Berlin 1959.

\bibit[2] D. Bao, S. S. Chern, and Z. Shen, {\it  An
Introduction to Riemann-Finsler Geometry,}  Springer, N.Y., Berlin 2000.

\bibit[3] L. Kozma and L. Tam{\' a}ssy,
Finsler geometry without line elements faced to applications,
{   \it Rep. Math. Phys.} {\bf 51} (2003), 233--250.

\bibit[4] L. Tam{\' a}ssy,
 Metrical almost linear connections in $TM$ for Randers spaces,
{\it Bull. Soc. Sci. Lett. Lodz Ser. Rech. Deform } {\bf 51} (2006), 147-152.

\bibit[5] Z. L. Szab{\' o},  All regular Landsberg metrics are Berwald,
{\it Ann Glob Anal Geom  } {\bf 34} (2008), 381-386.

\bibit[6] L. Tam{\' a}ssy,  Angle in Minkowski and Finsler spaces,
{\it Bull. Soc. Sci. Lett. Lodz Ser. Rech. Deform } {\bf 49} (2006), 7-14.

\bibit[7] G. S. Asanov,   Finsleroid-regular    space  developed.
   Berwald case, {\it  arXiv:} 0711.4180v1 [math.DG] (2007).

\bibit[8] G. S. Asanov,  Finsleroid-regular    space:~ curvature tensor,
continuation of gravitational Schwarzschild metric,
  {\it  arXiv}: 0712.0440v1 [math-ph] (2007).

\bibit[9] G. S. Asanov,
  Finsleroid-regular    space.  Gravitational  metric.  Berwald case,
{   \it Rep. Math. Phys.} {\bf 62} (2008), 103--128.

\bibit[10] G. S. Asanov,
 Finsleroid-regular   space.    Landsberg-to-Berwald implication,
{\it  arXiv:} 0801.4608v1 [math.DG],  (2008).

\end{document}